\documentclass{article} 
\usepackage{amssymb}
\usepackage{amsmath}
\begin{document} 
\input epsf.sty

\title{On polynomials orthogonal to all powers of a given
 polynomial on a segment}  
\author{F. Pakovich}
\date{}

\maketitle

%\tableofcontents
%\pagebreak

\def\be{\begin{equation}}
\def\ee{\end{equation}}
\def\bs{$\square$ \vskip 0.2cm}
\def\d{{\rm d}} 
\def\D{{\rm D}} 
\def\I{{\rm I}} 
\def\C{{\mathbb C}} 
\def\N{{\mathbb N}} 
\def\P{{\mathbb P}}
\def\Z{{\mathbb Z}}
\def\R{{\mathbb R}} 
\def\ord{{\rm ord}}

\def\e{\eqref}
\def\phi{{\varphi}}
\def\v{{\varepsilon}} 
\def\deg{{\rm deg\,}} 
\def\Det{{\rm Det}}
\def\dim{{\rm dim\,}} 
\def\Ker{{\rm Ker\,}} 
\def\Gal{{\rm Gal\,}}
\def\St{{\rm St\,}} 
\def\exp{{\rm exp\,}} 
\def\cos{{\rm cos\,}}
\def\circ{{\rm circ\,}} 
\def\diag{{\rm diag\,}} 
\def\GCD{{\rm GCD }}
\def\LCM{{\rm LCM }}
\def\mod{{\rm mod\ }}

\def\bp{\begin{proposition}}
\def\ep{\end{proposition}}
\def\bt{\begin{theorem}}
\def\et{\end{theorem}}
\def\be{\begin{equation}}
\def\l{\label}
\def\ee{\end{equation}}
\def\bl{\begin{lemma}}
\def\el{\end{lemma}}
\def\bc{\begin{corollary}}
\def\ec{\end{corollary}}
\def\pr{\noindent{\it Proof. }}
\def\note{\noindent{\bf Note. }}
\def\bd{\begin{definition}}
\def\ed{\end{definition}}

\newtheorem{theorem}{Theorem}[section]
\newtheorem{lemma}{Lemma}[section]
\newtheorem{definition}{Definition}[section]
\newtheorem{corollary}{Corollary}[section]
\newtheorem{proposition}{Proposition}[section]

\section{Introduction}

In this paper we investigate the following
``polynomial moment problem":
{\it for a complex polynomial $P(z)$ and distinct complex numbers $a,b$
to describe polynomials $q(z)$
such that}
\be \l{1}
\int^b_a P^i(z)q(z)\d z=0
\ee
{\it for all integer non-negative $i$.}
This problem attracted 
attention recently in the
series of papers \cite{bby}-\cite{by}, \cite{y1}, where (1) arose
as an infinitesimal version of the
center problem for the Abel differential equation with polynomial
coefficients in the complex domain.

The example of $P(z)=z$ shows that $q(z)=0$ can be the only
polynomial solution to
\eqref{1} since the Stone-Weierstrass
theorem implies that the only continuous complex-valued function
which is orthogonal to all powers of $z$ is zero. On the other
hand, if $P(a)=P(b),$ then non-trivial polynomial solutions to (1) always
exist. Indeed, it is enough to set
$q(z)=R^{\prime}(P(z))P^{\prime}(z),$ where $R(z)$ is any
complex polynomial.  Then for any $i\geq 0$ $$\int^b_a 
P^i(z)q(z)\d z=\int^{P(b)}_{P(a)}y^i R^{\prime}(y)\d y=0.$$

More generally, the following ``composition condition" imposed
on $P(z)$ and $Q(z)=\int q(z)\d z$ is sufficient for polynomials
$P(z),q(z)$ to satisfy (1): {\it there exist polynomials $\tilde
P(z),$ $\tilde Q(z),$ $W(z)$ such that}
\be \l{2}
P(z)=\tilde P(W(z)), \
\ Q(z)=\tilde Q(W(z)), \ {\it and} \ W(a)=W(b).
\ee
The
sufficiency of condition (2) follows from $W(a)=W(b)$ after the
change of variable $z \rightarrow W(z)$. It was suggested in the
papers \cite{bfy1}-\cite{bfy5} (``the composition conjecture") that,
under the additional assumption $P(a)=P(b),$
condition (1) is actually equivalent to   
condition (2). This conjecture is true in several special cases.
In particular, when $a,b$ are not critical points of $P(z)$   
(\cite{c}), when $P(z)$ is indecomposable (\cite{pa2}),
and in some other cases (see \cite{pry}, \cite{ro}, and the papers
cited above).

Nevertheless, in general the composition conjecture fails to be 
true; a class of counterexamples to the composition conjecture 
was constructed in \cite{pa1}. These counterexamples exploit
polynomials
which admit 
``double decompositions'' of the form
$P(z)=A(B(z))=C(D(z)),$ where $A(z),$ $B(z),$ $C(z),$ $D(z)$
are non-linear polynomials. If $P(z)$ is such a polynomial and,
in addition, the equalities $B(a)=B(b),$ $D(a)=D(b)$ hold, then
for any polynomial
$Q(z),$ which can be represented as $Q(z)=E(B(z))+F(D(z))$ for some
polynomials $E(z),F(z),$ condition (1) is satisfied with
$q(z)=Q^{\prime}(z)$ by linearity.
On the other hand, it can be shown (see \cite{pa1}) that if $\deg B(z)$ and 
$\deg D(z)$ are coprime then condition (2) is not satisfied already for
$Q(z)=B(z)+D(z)$. The simplest explicit counterexample
to the composition conjecture has
the following form:
$$ P(z)=T_6(z), \ \ \ \ q(z)=T_2^{\prime}(z)+T_3^{\prime}(z),\ \ \ \
a=-\sqrt{3}/2. \ \ \ \ b=\sqrt{3}/2,$$
where $T_n(z)$ denotes $n$-th Chebyshev polynomial.
%Note that, by the second Ritt theorem,
%double decompositions with $\deg A(z)=\deg D(z),$
%$\deg B(z)=\deg C(z)$ and $\deg B(z),\deg D(z)$
%coprime
%are equivalent either to decompositions with $$A(z)=z^nR^m(z),\ \ \ 
%B(z)=z^m, \ \ \ C(z)=z^m,\ \ \ D(z)= z^nR(z^m)$$ for a polynomial $R(z)$
%and ${\rm GCD}(n,m)=1$ or to decompositions with $$A(z)=T_m(z), \ \ \ B(z)=
%T_n(z),\ \ \ 
%C(z)= T_n(z), \ \ \ D(z)=T_m(z)$$
%for Chebyshev polynomials $T_n(z),$ $T_m(z)$ and ${\rm GCD}(n,m)=1$ (see
%\cite{ri},
%\cite{sch}).

The counterexamples above
suggest to transform the composition conjecture 
as follows (\cite{pa4}): {\it non-zero polynomials $P(z),$ $q(z)$ satisfy 
condition (1) 
if and only if $\int q(z)\d z$ can be represented as a {\it sum} of
polynomials $Q_j$ such that
\be \l{cc}
P(z)=\tilde P_j(W_j(z)), \ \ \
Q_j(z)=\tilde Q_j(W_j(z)), \ \ \ {\it and} \ \ \ W_j(a)=W_j(b)
\ee
for some $\tilde P_j(z), \tilde Q_j(z), W_j(z)\in \C[z]$.}
Note that we do not make any additional assumptions on the 
values of $P(z)$ at the points $a,b$ any more. 
In particular, the conjecture
implies that non-zero polynomials orthogonal to
all powers of a given polynomial $P(z)$ on $[a,b]$ exist if and only
if $P(a)=P(b).$ 
For $P(z)=T_n(z)$
conjecture \eqref{cc} was verified in \cite{pa3}.

This paper is organized as follows. We start from the description  of
necessary and sufficient
conditions for polynomials $P(z),q(z)$ to satisfy \eqref{1} in terms
of linear relations between
branches of the algebraic function $Q(P^{-1}(z)),$ where $P^{-1}(z)$ denotes
the algebraic function which is inverse to $P(z).$ More precisely,
let $P^{-1}_{i}(z),$ $1\leq i \leq n,$ be single-valued branches of $P^{-1}(z)$
in a simply-connected domain $U\subset \C$
containing no critical values of $P(z).$
We show that there exists a system of
equations 
\be \l{su}
\sum_{i=1}^nf_{s,i}Q(P^{-1}_{i}(z))=0, \ \ \ \ \ \ \ \ 1\leq s \leq \deg P(z),
\ee
where $f_{s,i}\in\{0,-1,1\},$
such that \eqref{1} holds if and only if \eqref{su} holds.
This system depends on $P(z),a,b$ and can be calculated explicitly
via a special graph $\lambda_P,$ called the ``cactus'' of $P(z),$
which contains
all the information about the monodromy of $P(z).$

The criterion
\eqref{su}
has a number of applications.
For example, it allows us
to reduce an infinite set of
equations \eqref{1} to a finite set of equations $v_k=0,$ ${0\leq k \leq M,}$
where $v_k$ are initial coefficients of
the Puiseux expansions of the combinations of branches in \eqref{su} and
$M$ depends only on degrees of $P(z)$ and $Q(z)$.
Furthermore, using the equivalence of \eqref{1} and \eqref{su},
we provide a variety of different conditions on a
collection $P(z),a,b$ 
under which conditions \eqref{1} and 
\eqref{2} are equivalent; 
such conditions are of interest because
of applications to the 
Abel equation (see \cite{bby}, \cite{bry}, \cite{by}). 
Essentially the finding of such conditions reduces to the finding conditions
under which system \eqref{su} implies that
\be \l{eb}
Q(P^{-1}_{i_1}(z))=Q(P^{-1}_{i_2}(z))
\ee 
for some $i_1\neq i_2.$
In its turn these conditions can be naturally given in terms
of combinatorics of the graph $\lambda_P.$

While an explicit form of system \eqref{su} depends on
the collection $P(z),a,b,$
there exists a necessary condition for \e{1} to be satisfied the form
of which is invariant with respect to $P(z),a,b.$
Namely, it is known (\cite{pa2}, \cite{pry}) that \e{1}
implies the equality  
\be \l{e1}
\frac{1}{d_{a}}\sum_{s=1}^{d_{a}}
Q(P_{a_s}^{-1}(z))=\frac{1}{d_{b}}\sum_{s=1}^{d_{b}}
Q(P_{b_s}^{-1}(z)),
\ee if $P(a)=P(b)$, or
the system 
\be \l{e2}
\frac{1}{d_{a}}\sum_{s=1}^{d_{a}}
Q(P_{a_s}^{-1}(z))=0, \ \ \ \ \ \ \ \ \ \ \frac{1}{d_{b}}\sum_{s=1}^{d_{b}}
Q(P_{b_s}^{-1}(z))=0,
\ee
if $P(a)\neq P(b),$ where
$P_{a_1}^{-1}(z),$
$P_{a_2}^{-1}(z), ... , P_{a_{f_a}}^{-1}(z)$
(resp. $P_{b_1}^{-1}(z),$ $P_{b_2}^{-1}(z),$ ... , \linebreak
$P_{b_{f_b}}^{-1}(z)$) denote
the branches
of $P^{-1}(z)$
which map points close
to $P(a)$ (resp. $P(b)$)
to points close to $a$ (resp. $b$).

In the fourth section of this paper we investigate equation \eqref{e1} and 
system \eqref{e2}. In particular, we establish a specific geometric
property of the monodromy groups of polynomials 
from which we deduce that if \e{e1} or \e{e2} is satisfied
for $P(z),Q(z)\in \C[z],$ $\deg P(z)=n, \deg Q(z)=m,$
then for coefficients of the Puiseux expansions near infinity
\be \l{ps2}
Q(P^{-1}_j(z))=\sum_{k=-m}^{\infty}
u_k\varepsilon_n^{jk}z^{-\frac{k}{n}}
\ee
the equality $u_k=0$ holds whenever $\GCD(k,n)=1.$
This fact agrees with conjecture \eqref{cc} and, in particular,
implies that for $P(z),q(z)$ satisfying \eqref{1} the numbers
$n$ and $m$ can not be coprime.
As an application of our analysis of \eqref{e1},\eqref{e2} 
we show in the fifth section that conditions \eqref{1} and \eqref{2}
are equivalent if at least one from points $a,b$ is not a critical
point of $P(z)$ or if $\deg P(z)=p^r$ for a prime number $p.$

Finally,
on the base of obtained results,
in the sixth section we show that for any
collection $P(z),a,b$ with $\deg P(z)<10$
conditions \eqref{1}
and \eqref{2} are equivalent except the case when 
$P(z),a,b$ is linearly equivalent to the collection
$T_6(z),- \sqrt{3}/2,\sqrt{3}/2.$ 
Since for $P(z)=T_n(z)$ all solutions to \eqref{1} were obtained
in \cite{pa3}, this provides
the complete solution of the polynomial moment problem for
$P(z),a,b$ with $\deg P(z)<10.$

\section{Criterion for a polynomial to be orthogonal to all
powers of a given polynomial} 

\subsection{Cauchy type integrals of algebraic functions} 
A quite general approach to the polynomial moment problem
was proposed in the paper \cite{pry} concerning Cauchy type integrals of
algebraic functions
\be \l{ci}
I(t)=I(\gamma,g,t)=\frac{1} {2\pi i}\int_{\gamma} \frac{g(z)dz}
{z-t}\,.
\ee
In this subsection we briefly recall it (see \cite{pry} for
details) and outline in this context the approach of this paper.

First of all notice that condition \eqref{1} is equivalent
to the condition
\be \l{1+}
\int^b_a P^i(z)Q(z)P^{\prime}(z)\d z=0
\ee
for $i\geq 0,$ where $Q(z)=\int q(z) \d z$ is normalized
by the condition ${Q(a)=Q(b)=0}$
($Q(a)$ always equals $Q(b)$
by \e{1} taken for $i=0$). 
Furthermore, vice versa, condition \e{1+}
with $Q(a)=Q(b)=0$ implies that \e{1} holds with $q(z)=Q^{\prime}(z).$
(Actually, it was the condition
\eqref{1+} that appeared initially in the papers
on differential equations cited above).

Indeed, condition \e{1+} is equivalent to the condition that
the function 
$$
H(t)= \int_a^b \frac{Q(z)P^{\prime}(z)dz}{P(z)-t} \ \ \ \ \
$$ vanishes identically near infinity, since near infinity 
$$H(t)=-\sum_{i=0}^{\infty}m_it^{-(i+1)}, \ \ \ \  {\rm where}
\ \ \ \  m_i=
\int^b_a P^i(z)Q(z)P^{\prime}(z)\d z.$$
On the other hand,
we have:  
\begin{multline} \l{mul}
\frac{dH(t)}{dt}=\int_a^b 
\frac{Q(z)P^{\prime}(z)dz}{(P(z)-t)^2}=
-\int_a^b Q(z)d\left(\frac{1}{P(z)-t}\right)=\\
=\frac{Q(a)}{P(a)-t}-\frac{Q(b)}{P(b)-t}+\tilde H(t),
\end{multline}
where
$$\tilde H(t)= \int_a^b \frac{q(z)dz}{P(z)-t}.$$ Since 
near infinity 
$$\tilde H(t)=-\sum_{i=0}^{\infty}\tilde m_it^{-(i+1)}, \ \ \ \  {\rm where}
\ \ \ \  \tilde m_i=
\int^b_a P^i(z)q(z)\d z,$$
it follows from \e{mul} that conditions \e{1} and \e{1+} are equivalent
whenever \linebreak
${Q(a)=Q(b)=0.}$

Furthermore, performing the change of variable
$z\rightarrow P(z),$
%the vanishing of $H(t)$
%near infinity becomes equivalent to
%the vanishing near infinity of Cauchy type integral
%\e{ci}, where $\gamma=P([a,b])$ and
%$g(z)$ is an algebraic function obtained by the analytic continuation
%of a germ of the algebraic function
%$g(z)=Q(P^{-1}(z))$ along $\gamma.$
we see that $H(t)$ coincides with integral \e{ci} where $\gamma=P([a,b])$ and
$g(z)$ is an algebraic function obtained by the analytic continuation
of a germ of the algebraic function
$g(z)=Q(P^{-1}(z))$ along $\gamma.$
Integral representation \eqref{ci}
defines a collection of univalent regular functions
$I_i(t);$ each $I_i(t)$ is defined in a domain $D_i$ of the complement
of $\gamma$ in $\C\P^1$. Denote by $I_{\infty}(t)$
the function defined in
the domain $D_{\infty}$ containing infinity. Then the vanishing of $H(t)$
near infinity becomes equivalent to the equality $I_{\infty}(t)\equiv 0.$

More generally,
consider integral \e{ci}, where $\gamma$ is a curve in the complex
plane $\C$ and $g(z)$ is any 
``piecewise-algebraic'' function on $\gamma$. More precisely, we assume that 
after removing from $\gamma$ a finite
set of points $\Sigma_{\gamma},$
the set $\gamma\setminus \Sigma_{\gamma}$ is an union
of topological segments $\cup \gamma_s$ such that 
for each $\gamma_s$ there exists a  
domain $V_s\supset \gamma_s$ and an analytic in $V_s$ 
algebraic function $g_s(z)$ such that $g(z)$
on $\gamma_s$ coincides with
$g_s(z).$ Furthermore, we 
assume that at the points of $\Sigma_{\gamma}$, complete analytic
continuations $\hat g_s(z)$ of $g_s(z)$ can ramify
but do not have 
poles. 
%Integral representation \eqref{ci}
%defines a collection of univalent regular functions
%$I_i(t);$ each $I_i(t)$ is defined in a domain $D_i$ of the complement
%of $\gamma$ in $\C\P^1$. Denote by $I_{\infty}(t)$
%the function defined in
%the domain $D_{\infty}$ containing infinity.
Below we sketch conditions for $I_{\infty}(t)$
to be a rational function; if these conditions are satisfied, then in
order to verify the equality 
$I_{\infty}(t)\equiv 0$
it is enough to examine possible poles.

Denote by $\Sigma_g$ the set 
of all singularities of $\hat g_s(z)$
in $\C\P^1.$
Show that any element $(I_i(t),U_i)$ can be analytically continued
along any curve
$S=S_{t_1,t_2}$ connecting points $t_1,t_2 \in \C\P^1$
and avoiding points from the sets $\Sigma_g$
and $\Sigma_{\gamma}.$ 
First of all notice that if $t_2\in \partial U_i$ then 
an analytical extension of $(I_i(t),U_i)$
to a domain containing $t_2$ 
is given simply by the integral
$I(\tilde\gamma,g,t),$ where $\tilde \gamma$ is a small
deformation of the $\gamma$
such that $t_2\in \tilde U_i.$
Furthermore, if $S=S_{t_1,t_2}$ is a simple curve connecting points
$t_1\in U_i,$ $t_2\in U_j,$ where
$U_i,$ $U_j$ are domains with a common segment of the boundary
$\gamma_{s}$ and $(g_s,V_s)$ is
the corresponding algebraic function, 
then the well-known 
boundary property of Cauchy type integrals (see e.g. \cite{mus}) implies that
$$(I_i(t),U_i\cap V_s)=(I_j(t),U_j\cap V_s)+(g_s,V_s).$$
Therefore, the analytic continuation of 
$(I_i(t),U_i)$ along $S$ can be defined via
the analytic continuation of the right side of this formula.

Finally,
for arbitrary domains $U_i,$ $U_j$ and
a curve $S=S_{t_1,t_2}$ connecting points $t_1\in U_i$ and $t_2\in U_j,$ 
the analytic continuation $(I_i(t),U_i)_S$ of 
the element $(I_i(t),U_i)$ along $S$ can be defined inductively as follows.
Let $S\cap \gamma =\{c_1,c_2, ... ,c_r\},$ 
$c_s\in V_s,$ $1\leq s \leq r,$
and 
let $(g_1,V_1), (g_2,V_2), ... , (g_r,V_r)$ be 
the corresponding algebraic functions.
Define a germ $g_{\gamma,S}$ of an algebraic function near the point 
$t_2$ by the formula:
$$g_{\gamma,S}=\sum_{i=1}^r (g_i,V_i)_{S_{c_i}},  $$
where $(g_i,V_i)_{S_{c_i}}$, $1\leq i \leq r,$ denotes the analytic 
continuation
of the element $(g_i,V_i)$ (taken with the sign corresponding
to the orientation of the crossing of $S$ and $\gamma$)
along a part of $S$ from $c_i$ to $t_2.$ 
Then by induction we have:
\be \l{skl}
(I_i(t),U_i)_{S}=(I_j(t),U_j)+g_{\gamma,S}.
\ee 
In particular, a complete analytic continuation $\hat I_i(t)$ of the element
$(I_i(t),U_i)$ is
a multi-valued analytic function with a finite set of
singularities $\Sigma_{\hat I_i} \subset \Sigma_g \cup \Sigma_{\gamma}.$

From formula \e{skl} one deduces the following criterion
(\cite{pry}):
{\it $\hat I_i(t)$ is a rational function 
if and only if the equality
\be \l{us}
g_{\gamma,S}=0
\ee
holds for any curve
$S=S_{t_1,t_2}$ as above with $t_1=t_2\in U_i.$}
Indeed, the necessity of 
\e{us} is obvious. To establish the sufficiency
observe that \e{us} implies, in particular, that $\hat I_i(t)$ has
no ramification in its singularities. 
Since
a Cauchy type integral near the end of the line of integration $u$
(and, similarly, near points of $\Sigma_{\gamma}$) is
a sum of a part which has a logarithmic branching at $u$
with a part that bounded at
$u$ (see e.g. \cite{mus}) and
$\hat g_s(z)$ do not have 
poles at $\Sigma_{\gamma},$
this fact implies that the singularities of $\hat I_i(t)$
located on $\gamma$ are removable.
On the other hand, if a singularity
$t_0$ of $\hat I_i(t)$ is contained in $\C\P^1\setminus \gamma,$ then
formula \e{skl} implies that
$t_0\in \Sigma_g.$ Since $\hat I_i(t)$ has no ramification
at $t_0$ it follows that in this case $t_0$ 
is a pole the worst.

Although the method above in principle is constructive
its practical application is rather
difficult since the calculation of sums $g_{\gamma,S}$ 
is complicated.
In this paper we propose a modification
of the method above
designed
specially for the polynomial moment problem.
This modification
permits to avoid any analytic continuations  
and allows us to obtain a necessary and sufficient conditions
for equality \e{1} to be satisfied in 
a closed and convenient form.
The idea is
to choose a very special way of integration $\Gamma$ connecting points
$a,b$ (we can use any of them since integrals \e{1}
do not depend on $\Gamma$). It turns out that 
$\Gamma$ can by chosen  
so that $\C\P^1\setminus P(\Gamma)$ consists of
a {\it unique} domain. Then condition $I_{\infty}(t)\equiv 0$
simply reduces to the condition that
the corresponding algebraic functions $g_s(z)$ vanish on $P(\Gamma).$
Furthermore, we choose $\Gamma$ as a subset of a special tree $\lambda_P,$
called the cactus of $P(z),$ which contains
all the information about the monodromy $P(z).$ 
This allows us to relate properties of the collection
$P(z),a,b$ which are connected to the polynomial moment 
problem to combinatorial properties of the pair consisting of
the tree $\lambda_P$ and the   
path 
connecting points $a,b$ on $\lambda_P.$

\subsection{Cacti} \l{cacti}
To visualize the monodromy group of a polynomial $P(z)$
it is convenient to consider a graphical object $\lambda_P$ called the
{\it cactus} of $P(z)$ (see e. g. \cite{cac1}, \cite{cac2}). 

Let $c_1,c_2, ... ,c_k$ be all
finite critical values of $P(z)$ and let $c$ be a not critical value. Draw 
a star $S$ joining $c$ with $c_1,c_2, ... ,c_k$ by non intersecting arcs
$\gamma_1, \gamma_2, ... ,\gamma_k.$ We will suppose that $c_1,c_2, ...
,c_k$ are numerated in such a way that in a counterclockwise rotation
around $c$ the arc $\gamma_s$, $1 \leq s \leq k-1,$ is followed by the arc
$\gamma_{s+1}.$ By definition, the cactus $\lambda_P$ is the preimage of $S$
under the map $P(z)\,:\, \C \rightarrow \C.$ More precisely, we consider
$\lambda_P$ as a $k+1$-colored graph embedded into the Riemann sphere:
vertices of $\lambda_P$ colored by the $s$-th color, where $1\leq s \leq 
k,$ are preimages of the point $c_s,$ vertices colored by the $k+1$-th 
color (to 
be definite we will suppose that it is the white color) are 
preimages
of the point $c,$ and edges of $\lambda_P$ are preimages of the arcs
$\gamma_s,$ $1 \leq s \leq k.$ It is not difficult to show that the graph
$\lambda_P$ is connected and has no cycles. Therefore,
$\lambda_P$ is a plane tree.  

The valency of a non-white vertex $z$ of $\lambda_P$
coincides with the multiplicity of $z$ with respect to $P(z)$ while
all white vertices of $\lambda_P$ are of the same valency
$n=\deg P(z).$ The set of all edges of $\lambda_P$ adjacent to a
white vertex $w$ is called by a {\it star} of $\lambda_P$ centered at $w.$
Clearly, $\lambda_P$ has $nk$ edges and $n$ stars. The set of stars
of $\lambda_P$ is naturally identified with the set of branches of $P^{-1}(z)$
as follows.
Let $U$ be a
simply connected domain containing no critical values of $P(z)$ such 
that $S\setminus\{c_1,c_2,
... ,c_k\}\subset U.$ By the monodromy theorem in $U$ there exist
$n$ single valued branches
of $P^{-1}(z).$
Any such a branch $P_i^{-1}(z),$ $1\leq i \leq n,$
maps $S\setminus\{c_1,c_2, ... ,c_k\}$ into a star of
$\lambda_P$ and we will label the corresponding star by the symbol $S_i$
(see Fig. 1).

\begin{figure}
\medskip
\epsfxsize=12truecm
\centerline{\epsffile{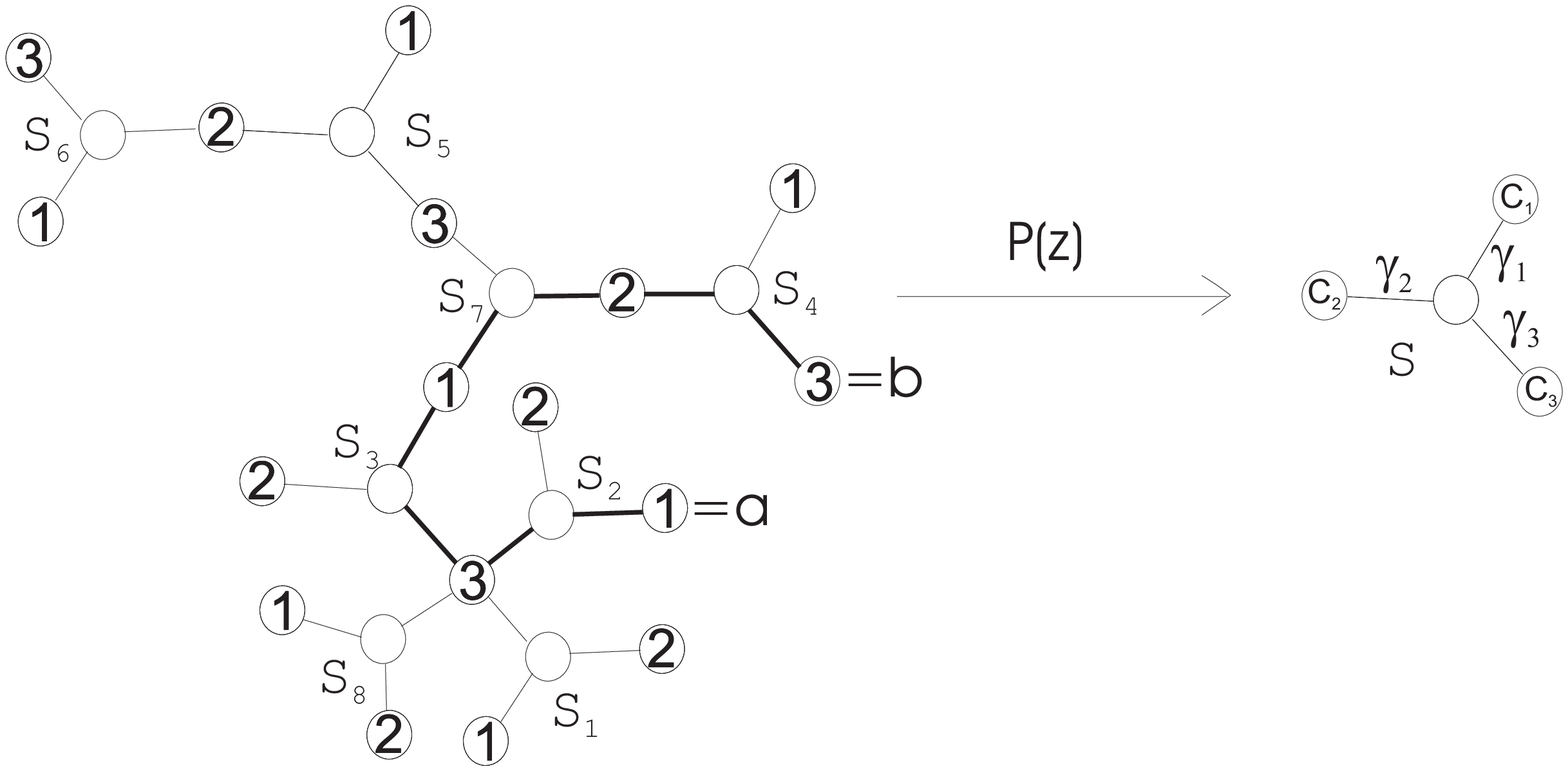}}
\smallskip
\centerline{Figure 1}
\medskip
\end{figure}

The cactus $\lambda_P$ permits to reconstruct the monodromy group
$G_P$ of $P(z).$ Indeed, $G_P$ is generated by  
the permutations $g_{s}\in S_n,$ $1\leq s \leq k,$ where 
$g_{s}$ is defined by the condition
that the analytic continuation of
the element $(P_i^{-1}(z),U),$ $1\leq i \leq n,$ along a counterclockwise
oriented loop $l_s$ around $c_s$
is the element $(P_{g_s(i)}^{-1}(z),U).$ Having in mind the identification of 
the set of stars of $\lambda_P$ with the set of branches of $P^{-1}(z),$
the permutation $g_s,$ $1\leq s \leq k,$ can be 
identified with the permutation $\hat g_s,$
$1\leq s \leq k,$
acting
on the set of starts of $\lambda_P$ in the following way: $\hat g_s$ sends the
star $S_i,$ $1\leq i \leq n,$ to the ``next'' one in a counterclockwise
direction around its vertex of color $s.$ For example,
for the cactus shown on Fig. 1 we have: $g_1=(1)(2)(37)(4)(5)(6)(8),$
$g_2=(1)(2)(3)(47)(56)(8),$ $g_3=(1238)(4)(57)(6).$

Note that since $P(z)$ is a polynomial, the permutation
$g_{\infty}=g_1g_2 ... g_k$ is a cycle of
length $n.$ Usually, we will choose the numeration of $S_i,$ $1\leq i \leq
n,$ in such a way that this cycle coincides with the cycle $(12...n).$

\subsection{Criterion}\l{criter} In this subsection
we give explicit necessary and sufficient conditions for $P(z),$ $q(z)\in
\C[z]$ and $a,b\in \C,$ $a\neq b,$ to satisfy \e{1}, \e{1+}.
For this propose we choose the way of integration $\Gamma_{a,b}$ connecting
$a,b$ so
that $\Gamma_{a,b}$ would be a subset of $\lambda_P.$

More precisely, for any $P(z)\in \C[z]$ and $a,b\in \C$ let us define
an {\it extended} cactus
$\tilde \lambda_P=\tilde \lambda_P(c_1,c_2,\, ... \,,c_{\tilde k})$
as follows. Let $c_1,c_2, ... ,c_{\tilde k}$ be all finite critical values
of $P(z)$ complemented by $P(a)$ or $P(b)$ (or by both of them) if $P(a)$
or $P(b)$ is not a critical value of $P(z).$ Consider an extended star
$\tilde S$ connecting $c$ with $c_1,c_2, ... ,c_{\tilde k}$ and set
$\tilde \lambda_P=P^{-1}\{\tilde S\}$ (we suppose that $c$ is chosen
distinct from $P(a),P(b)$). Clearly, $\tilde \lambda_P$ considered as a
$\tilde k+1$ colored graph is still connected and has no cycles.
Furthermore,
by construction the points $a,b$ are vertices of $\tilde \lambda_P.$ Since
$\tilde \lambda_P$ is connected there exists an oriented path $\Gamma_{a,b}
\subset \tilde \lambda_P$
with the starting point $a$ and the ending point $b.$ Moreover, since
$\tilde \lambda_P$ has no cycles there exists exactly one such a path. We
choose $\Gamma_{a,b}$ as a new way of integration.

Let $U$ be a domain as above and let $Q(z)=\int q(z) \d z$ be normalized
by the condition ${Q(a)=Q(b)=0}.$ 
For each
$s,$ $1\leq s \leq \tilde k,$ define a linear combination
$\phi_s(z)$ of branches $Q(P^{-1}_{i}(z)),$ $1\leq i \leq n,$ in $U$
as follows. Set 
\be \l{piz}
\phi_s(z)=\sum_{i=1}^nf_{s,i}Q(P^{-1}_{i}(z)), 
\ee
where $f_{s,i}\neq 0$ if and only if the path
$\Gamma_{a,b}$ passes through a vertex $v$ of the star
$S_i$
colored by the
$s$-th color
(we do not take into account the stars $S_i$ for which $\Gamma_{a,b}\cap S_i$
contains only the point $v$). Furthermore, if under a moving along $\Gamma_{a,b}$ the vertex $v$
is followed by the center of $S_i$ then $f_{s,i}=-1$ otherwise
$f_{s,i}=1$.  As an example consider the cactus shown on Fig. 1. Then 
for the path
$\Gamma_{a,b}$ pictured by the fat 
line we 
have:  $$\phi_1(z)=
-Q(P^{-1}_{2}(z))+Q(P^{-1}_{3}(z))- Q(P^{-1}_{7}(z)),$$ $$\phi_2(z)=
Q(P^{-1}_{7}(z))-Q(P^{-1}_{4}(z)),$$ $$\phi_3(z)=
Q(P^{-1}_{2}(z))-Q(P^{-1}_{3}(z))+ Q(P^{-1}_{4}(z)).$$
\vskip 0.1cm
\bt \l{t1}
Let $P(z),q(z)\in \C[z]$, $a,b\in \C,$ $a\neq b,$
and let $\tilde \lambda_P(c_1,c_2,\, ... \,,c_{\tilde k})$ be an extended
cactus corresponding to the collection $P(z),a,b.$
Then 
\eqref{1} holds if and only if 
the equality $\phi_s(z)\equiv 0$ holds
in $U$ for any $s,$ $1\leq s \leq \tilde k.$
\et 

\pr Indeed,
condition \eqref{1+} is equivalent to the condition that the function
$$
H(t)= \int_{\Gamma_{a,b}} \frac{Q(z)P^{\prime}(z)dz}{P(z)-t} \ \ \ \ \
$$ vanishes identically near infinity.
On the other hand, using
the change of variable $z\rightarrow P(z),$ we can express the function
$H(t)$ as a sum of Cauchy type integrals of algebraic functions as follows:
\be \l{f}
H(t)=\sum_{s=1}^{\tilde k}\int_{\gamma_s}\frac{\phi_s(z)}{z-t}\, \d
z.
\ee
Since this formula implies that $H(t)$ is analytic in a domain
$V=\C\P^1\setminus S$ we see that the vanishing of $H(t)$
near infinity is equivalent to the
condition that $H(t)\equiv 0$ in $V.$

Let $z_0$ be an interior point of $\gamma_s,$ $1\leq s \leq \tilde k.$ 
Then by the well-known boundary property of
Cauchy type integrals (see e.g. \cite{mus}) we have:

$$
\lim_{t \to z_0}\!\!\!\!\,^+H(t)-\lim_{t \to z_0}
\!\!\!\!\,^-H(t)=\phi_s(t_0),
$$
where the limits are taken respectively for $t$ 
tending to $z_0$ from the ``left'' and 
from the ``right'' parts of $V$ with respect to $\gamma_s.$ 
If $H(t)\equiv 0$ in $V,$ then 
$$
\lim_{t \to z_0}\!\!\!\!\,^+H(t)=\lim_{t \to z_0}
\!\!\!\!\,^-H(t)=0,
$$
and, therefore, $\phi_s(z_0)= 0.$
Since this equality holds for any
interior point $z_0$ of any arc $\gamma_s,$ $1\leq s \leq \tilde k,$ 
we conclude that  
$\phi_s(z)\equiv 0,$ $1\leq s \leq \tilde k,$ in $U.$ On the other hand,
if $\phi_s(z)\equiv 0,$ $1\leq s \leq \tilde k,$ in $U,$ then
it follows directly from formula \e{f} that
$H(t)\equiv 0$ in $V.$

\vskip 0.2cm
Note that some of equations $\phi_s(z)\equiv 0,$ $1\leq s \leq \tilde k,$
could be trivial. 
This happens exactly for the values $s$ such that
the path $\Gamma_{a,b}$ does not pass through vertices
colored by the $s$-th color.
Note also that
equations \e{piz}
are linearly dependent. Indeed, 
for each $i$ such that there exists an index $s,$ 
$1\leq s
\leq \tilde k,$ with $f_{s,i}\neq 0$ there exist exactly two such
indices $s_1,s_2$ and $c_{s_1,i}=-c_{s_2,i}.$ Therefore, the equality
$$\sum_{s=1}^{\tilde k}\phi_s(t)=0$$
holds in $U$.

\subsection{Checking the criterion}
Let $P(z),Q(z)\in \C[z],$ $\deg P(z)=n,$ $\deg Q(z)=m.$
Let $U$ be a simply connected domain containing no critical values of $P(z)$
and let $P_i^{-1}(z),$ $1\leq i \leq n,$ be branches of $P^{-1}(z)$
in $U.$
In this subsection we provide a simple estimation 
for the order of a zero in $U$ of a function of the form
$$\psi(z)=\sum_{i=1}^nf_iQ(P^{-1}_i(z)), \ \ \ \ \ \ \ \ \ \ \ \ f_i\in \C,$$
via the degrees of $P(z)$ and $Q(z).$ This reduces 
the verification of the criterion
to the calculation of a finite set of initial coefficients of
Puiseux expansions of functions \e{piz}
and, as a corollary, provides a practical method for checking an infinite
set of equation \e{1} in a finite number of steps.

\bl\l{ur} If $\psi(z)\neq 0$ then $\psi(z)$
satisfies an equation
\be \l{urr}
y^N(z)+a_1(z)y^{N-1}(z)+\ ... \ + a_N(z)=0,
\ee
where $a_j(z)\in \C[z],$
$a_N(z)\neq 0,$ and $N\leq n!.$
Furthermore, 
$\deg a_j(z)\leq \left(\frac{m}{n}\right)^{j},$ $1\leq j \leq N.$
\el

\pr Indeed, if $\psi(z)\neq 0$ then, since $\psi(z)$ is a sum of algebraic
functions, $\psi(z)$ itself is an algebraic function and
therefore satisfies an algebraic equation \e{urr}
with $a_i(z)\in \C(z),$ $1\leq j \leq N.$ Furthermore, we can
suppose that this equation
is irreducible. Then $a_N(z)\neq 0$
and the number $N$ 
coincides with the number of different analytic continuations
$\psi_j(z)$ of
$\psi(z)$ along closed curves. Clearly,
$N$ can be bounded by the number $N_1$ of different elements of
the monodromy group of $P(z).$ In its turn, $N_1$ is bounded by the
number of elements of the full symmetric group $S_n.$ Hence, $N\leq n!.$

Furthermore, since $P(z),Q(z)$ are polynomials, the rational functions
$a_j(z),$ $1\leq j \leq N,$ as the symmetric functions
of $\psi_j(z),$ $1\leq j \leq N,$ have no poles in $\C$ and therefore
are polynomials. Finally, since near infinity branches $P^{-1}_i(z),$
$1\leq i \leq n,$ 
of $P^{-1}(z)$
are represented 
by the Puiseux series 
\be 
P^{-1}_i(z)=\sum_{k=-1}^{\infty}v_k\varepsilon_n^{ik}z^{-\frac{k}{n}},
\ \ \ \ \ \ \ \ v_k\in \C, \ \ \ \ \ \ \ \ \varepsilon_n=exp(2\pi i/n),
\ee
the first non-zero exponent of the Puiseux series 
at infinity for the functions
$\psi_j(z),$ $1\leq j \leq N,$ 
is less or equal than $m/n.$ It follows that
$\deg a_j(z)\leq \left(\frac{m}{n}\right)^{j},$ $1\leq j \leq N.$

\bc Let $z_0 \in U.$ To verify that $\psi(z) \equiv 0$ it is enough
to check that the first $\left(\frac{m}{n}\right)^{n!}+1$ coefficients of
the series
$\psi(z)=\sum_{k=0}^{\infty} v_k (z-z_0)^k$ vanish.

\ec

\pr Indeed, suppose that $\ord_{z_0}\psi(z) > \left(\frac{m}{n}\right)^{n!}$
but $\psi(z)\neq 0.$ Then,
by lemma \ref{ur}, $\psi(z)$ satisfies \e{urr}, where 
$\deg a_j(z)\leq \left(\frac{m}{n}\right)^{j}\leq \left(\frac{m}{n}
\right)^{n!},$ $1\leq j \leq N,$ and $a_N\neq 0.$
It follows that
$$\ord_{z_0}\{\psi^N(z)\} >
\ord_{z_0}\{a_{i_1}(z)\psi^{N-i_1}(z)\} >\ ... \ >\ord_{z_0}\{a_{i_k}(z)
\psi^{N-i_k}(z)\},$$ where $0 < i_1 < i_2 <\ ... \ <i_k=N$ are all indices
for which $a_i(z) \neq 0.$
Therefore,
$$\ord_{z_0}\{\psi^N(z)+
a_1(z)\psi^{N-1}(z) +\ ... \ +a_{N}(z)\}
=\ord_{z_0}\{a_{N}(z)\}<\infty $$ in contradiction with equality
\e{urr}.

\section{Application to a description of definite polynomials}
In this section, as a first
application of theorem \ref{t1}, we provide
a number of conditions
on a collection $P(z),a,b,$ where $P(z)\in \C[z],$ $a,b\in \C,$
$a\neq b,$ under which conditions \e{1} and \e{2}
are equivalent; such collections are called {\it definite} and
are of interest because of applications to the Abel equation (see \cite{bby},
\cite{bry}, \cite{by}).

\subsection{A combinatorial condition for a change of variable}

The simplest form of the equality 
$\phi_s(z)=0$ 
is equality \e{eb}.
Furthermore, \e{eb} has a clear compositional meaning.

\bl \l{lcomp} The equality \eqref{eb} holds if and only if
\be \l{comp}
P(z)=\tilde P(W(z)), \ \ \ \ \ \ \ \ Q(z)=\tilde Q(W(z))
\ee
for
some polynomials $\tilde P(z),$ $\tilde Q(z),$ $W(z)$ with $\deg W(z)>1.$
\el

The proof of this lemma easily follows from the L\"{u}roth theorem
(see e.g. \cite{pa2}, \cite{ro}).
If condition \e{comp} is satisfied we say that
{\it polynomials $P(z),$ $Q(z)$ have a (non-trivial)
common right divisor in the composition algebra of polynomials.}

Below we give a convenient combinatorial
condition on a collection $P(z),a,b$ which implies that for any 
$q(z)$ satisfying \e{1} polynomials $P(z),$ $Q(z)=\int q(z)\d z$ have a
common right divisor in the composition algebra of polynomials.
The use of this condition permits, after the change of variable
$z\rightarrow W(z),$
to reduce the solution of the polynomial
moment problem for a polynomial $P(z)$ to that for a polynomial
of lesser degree $\tilde P(z).$

Let $\tilde \lambda_P$ be a $\tilde k+1$ colored extended cactus
corresponding to a collection $P(z),a,b$ and let
$\Gamma_{a,b}\subset \lambda_P$
be the path connecting points $a,b$ on $\tilde \lambda_P.$ 
For each $s,$ $1\leq s \leq \tilde k,$
define the weight $w(s)$ of the $s$-th color
on $\Gamma_{a,b}$ as a number of vertices $v\in \Gamma_{a,b}$
colored by the $s$-th color with the convention that vertices $a,b$ are
counted with the coefficient $1/2.$ For example, for $\Gamma_{a,b}$
shown on Fig. 1 we have $w(1)=w(3)=3/2,$ $w(2)=1.$ 

\bt \l{rt}  Let $P(z),q(z)\in \C[z]$, $q(z)\neq 0,$
$a,b\in \C,$ $a\neq b$ satisfy \e{1}. Suppose that there exists
$s,$ $1\leq s \leq \tilde k,$ such that $w(s)=1$ on $\Gamma_{a,b}.$
Then $P(z),Q(z)$ have a common right divisor in the composition algebra.
\et
\pr Indeed, the construction of $\Gamma_{a,b}$ implies that if
$w(s)=1,$ then $f_{s,i}\neq 0$
exactly for two values $i_1,i_2,$ $1\leq i_1,i_2 \leq n.$ Moreover, for these values we have
$c_{s,{i_1}}=-c_{s,{i_2}}$ and
hence the equality $\phi_s(z)=0$
reduces to \e{eb}. Therefore,
$P(z)$ and $Q(z)$ have a common right divisor in the composition algebra
by lemma \ref{lcomp}.

\subsection{Reduction}

Although condition \eqref{comp} in general is weaker than condition \eqref{2} 
it turns out that in order to prove
that for any collection
$P(z),a,b,$ $a\neq b,$
satisfying some condition $\mathcal R$
conditions \eqref{1} and \eqref{2} are equivalent
it is often enough to show 
that for any such a collection
condition \e{1} implies condition \e{comp}.
Say that a condition $\mathcal R$
is {\it compositionally stable} if for any
collection $P(z),a,b,$ $a\neq b,$ satisfying $\mathcal R$
such that $P(z)=\tilde P(W(z))$ for
some $\tilde P(z), W(z)\in \C[z],$ ${\deg W(z) >1,}$ $W(a)\neq W(b),$   
the collection $\tilde P(z),W(a),$ $W(b)$ also
satisfies $\mathcal R.$ For instance, the following condition is
compositional stable: at least one point from $a,b$ is
not a critical point of $P(z).$ An other example of a
compositional stable condition
is the following one: $\deg P(z)=p^r,$ where $p$ is a prime.

\vskip 0.2cm

\bl \l{rl}
Let $\mathcal R$ be
a compositionally stable condition.
Suppose that for any collection $P(z),a,b,$ $a\neq b,$
satisfying $\mathcal R,$
condition \e{1} implies condition \e{comp}.
Then for any collection $P(z),a,b,$ $a\neq b,$
satisfying $\mathcal R$ conditions \e{1} and \e{2} are equivalent.
\el

\vskip 0.2cm

\pr
Let $P(z),a,b$ be a collection satisfying $\mathcal R.$
Suppose that \e{1} holds for some $q(z)\in \C[z].$ Then it follows from
\e{comp} that $\C(P,Q)$ is a proper
subfield of $\C(z).$ Therefore, by the L\"{u}roth theorem
\be\l{lu}
\C(P,Q)=\C(W_1)
\ee for some rational function $W_1(z),$ $\deg W_1(z) > 1,$ 
and without loss of generality we can assume that
$W_1(z)$ is a polynomial. To prove the lemma it is enough to show
that the equality $W_1(a)=W_1(b)$ holds. 

Observe that equality \e{lu} is equivalent to the statement
that 
\be \l{comp+}
P(z)=P_1(W_1(z)), \ \ \ \ \ \ \ \ Q(z)=Q_1(W_1(z))
\ee
for
some polynomials $P_1(z),$ $Q_1(z)$ such that
$P_1(z)$ and $Q_1(z)$ have no a
common right divisor in the composition algebra.
On the other hand, 
performing the change of variable $z\rightarrow W_1(z)$
we see that condition \e{1} is satisfied also for
$P_1(z), Q_1^{\prime}(z), W_1(a),W_1(b).$
Therefore, if $W_1(a)\neq W_1(b),$ then,
taking into account that $\mathcal{R}$ is compositionally stable,
we find that $$P_1(z)=P_2
(W_2(z)), \ \ \ \ \ \ \ Q_1(z)=
Q_2(W_2(z))$$ for some $P_2(z), Q_2(z), W_2(z)\in \C[z]$ with
$\deg W_2(z) >1$ that contradicts the fact that
$P_1(z),$ $Q_1(z)$ have no a
common right divisor in the composition algebra.

\vskip 0.2cm

\vskip 0.2cm

\subsection{Description of some classes of definite
polynomials}

As a first application of theorem 
\ref{rt} and lemma \ref{rl} we give a simple proof
of the following assertion conjectured in \cite{pry}.

\bc \l{c1} Let $P(z),q(z)\in \C[z]$, $q(z)\neq 0,$
$a,b\in \C,$ $a\neq b.$ Suppose that
$P(a)=P(b)=c_1$ and that all the
points of the preimage
$P^{-1}(c_1)$ except possibly $a,b$ are not critical points
of $P(z).$ Then conditions \e{1} and \e{2} are equivalent. 
\ec

\pr
Since the chain rule implies that
the condition of the corollary is compositionally
stable it is enough to show that 
$P(z),$ $Q(z)$ have a
common right divisor in the composition algebra.
To establish it observe that $\Gamma_{a,b}$ can not pass through
vertices of $\lambda_P$ of the valency $1$ distinct from $a,b.$
Therefore, the
condition of the corollary implies 
that $w(1)=1.$ It follows now from theorem \ref{rt} that
$P(z),$ $Q(z)$ have a
common right divisor in the composition algebra.

\vskip 0.2cm A slight modification of the idea used in the proof
of corollary \e{c1}
leads to the following
statement.

\bc Let $P(z),q(z)\in \C[z]$, $q(z)\neq 0,$
$a,b\in \C,$ $a\neq b.$ Suppose that
$P(a)=P(b)=c_1$ and that for any critical
value $c$ of $P(z)$ except possibly $c_1$ the preimage $P^{-1}(c)$
contains only one critical point. Then conditions \e{1} and \e{2} are 
equivalent. \ec

\pr
Again, it follows from the chain rule that 
the condition of the corollary is compositionally stable.
Furthermore, observe that
the path $\Gamma_{a,b}$ contains at least one vertex $v$ 
of a color $s\neq 1.$  
Since $\Gamma_{a,b}$ can not pass through vertices of the
valency 
$1$ distinct from $a,b,$ it follows from the condition
of the corollary that 
the equality
$w(s)=1$ holds and, therefore, by theorem \ref{rt},
$P(z),$ $Q(z)$ have a
common right divisor in the composition algebra.

\vskip 0.2cm

Finally, we give a new proof of an
assertion from the paper \cite{pry}
which provides some geometric condition for a collection $P(z),a,b$
to be definite. It turns out that
this assertion actually also can be regarded as a particular
case of theorem \ref{rt}.
For a curve $M$ 
denote by $V_{M,\infty}$
the domain from the collection of domains $\C\P^1\setminus M$ which contains infinity. 
For an oriented curve $L$ and points 
$d_1,d_2\in L$ 
denote by $L_{d_1,d_2}$ the part of
$L$ between $d_1$ and $d_2.$ 

\bc Let $P(z),q(z)\in \C[z]$, $q(z)\neq 0,$ 
$a,b\in \C,$ $a\neq b.$ Suppose that
$P(a)=P(b)=c_1$ and that
there exists a curve $L$ 
connecting points $a,b$
such that $c_1$
is a simple point of $P(L)$ and $c_1\in \partial V_{P(L),\infty}.$
Then conditions \e{1} and \e{2} are 
equivalent.
\ec

\pr
We will keep the notation introduced in subsection \ref{cacti}
and \ref{criter}.
Let $a^+$ (resp. $b^-$) be a point on $L$ near the point $a$
(resp. $b$) and
let $U$ be a
simply connected domain containing no critical values of $P(z)$ such 
that the sets $S\setminus\{c_1,c_2,... ,c_{\tilde k}\},$
$P(L_{a,a^+})\setminus c_1,$
and $P(L_{b^-,b})\setminus c_1$ are subsets of
$U.$ Recall that there is a natural correspondence between branches
$P_i^{-1}(z),$ $1\leq i \leq n,$ of
$P^{-1}(z)$ in $U$ and stars of the cactus $\lambda_P:$
branch $P_i^{-1}(z)$ maps $U$ on a domain $U_i$ containing $S_i.$
Denote by $U_{j_1}$ (resp. $U_{j_2}$) the domain containing the point $a^+$
(resp. $b^-$). Then by construction the result of the analytic continuation
of the element $(P^{-1}_{j_1}(z),U)$ along the curve
$P(L_{a^+,b^-})$ is the element $(P^{-1}_{j_2}(z),U).$
\begin{figure}
\medskip
\epsfxsize=4.5truecm
\centerline{\epsffile{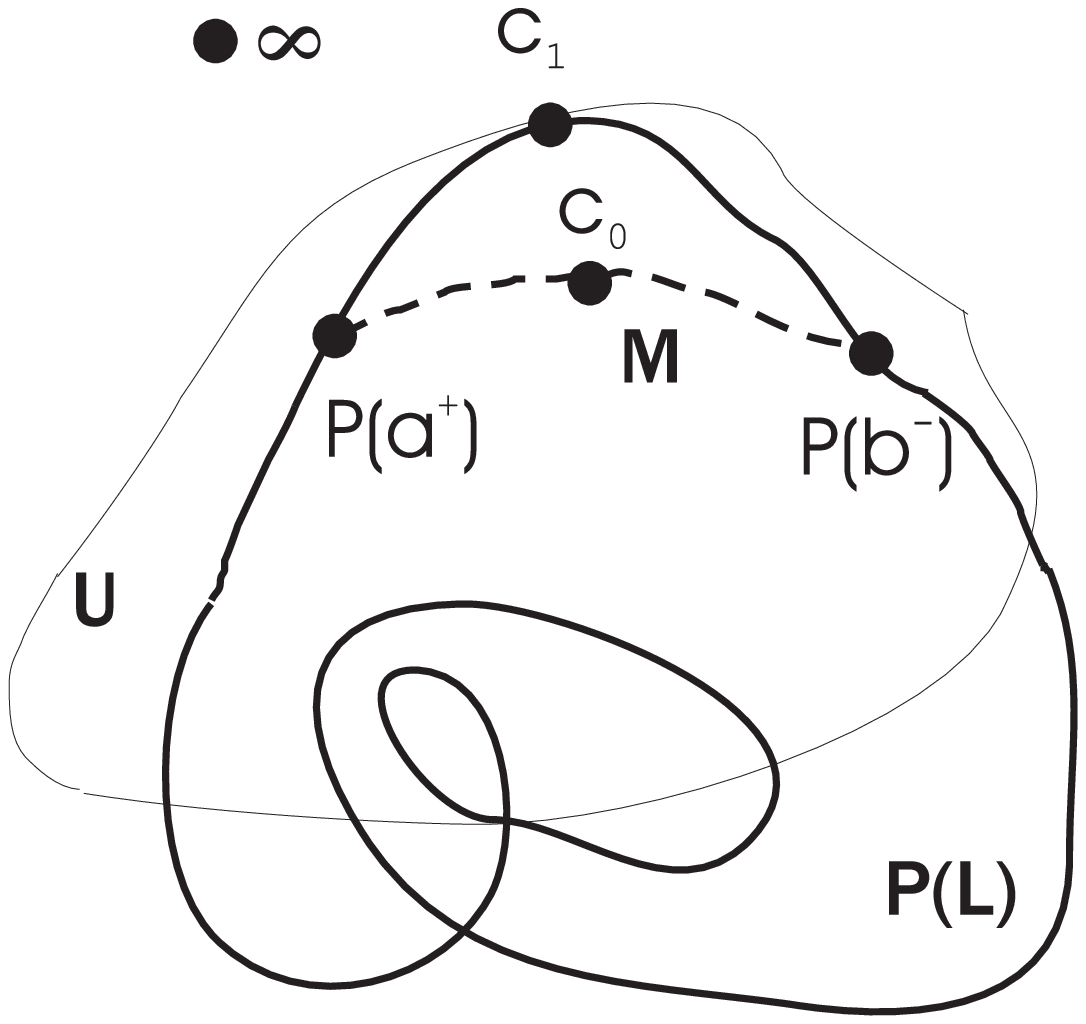}}
\smallskip
\centerline{Figure 2}
\medskip
\end{figure}
Let $c_0$ be an interior point of $U$ close to $c_1.$
Consider a small deformation $M$ of the curve $P(L_{a,b})$ obtained
as follows: change the part of $P(L_{a,b})$ connecting $c_1$
and $P(a^+)$
to an arc $\gamma^+\subset U$ connecting $c_0$ with $P(a^+)$
and, similarly, change the part of $P(L_{a,b})$ connecting
$P(b^-)$ and $c_1$ 
to an arc $\gamma^-\subset U$ connecting $P(b^-)$ and $c_0$ (see Fig. 2).

Let now $l_M=l_{i_1}^{j_1}l_{i_2}^{j_2}\, ...\, l_{i_r}^{j_r}$ 
be the image of the curve $M$ in 
$\pi_1(X, c_0),$ where $X=\C\P^1\setminus \{c_1,c_2,c_3, \, ... \,
c_{k}\}.$ 
Since the result of the analytic continuation
of the element $(P^{-1}_{j_1}(z),U)$ along the curve
$M$ is still the element $(P^{-1}_{j_2}(z),U),$
the final element of the chain of stars 
$$\Omega=<S_{j_1}\,,\;S_{g_{i_1}^{j_1}(j)}\,,\;
S_{g_{i_1}^{j_1}g_{i_2}^{j_2}(j)}
\,,\;
S_{g_{i_1}^{j_1}g_{i_2}^{j_2}g_{i_3}^{j_3}(j)}\,,\;...\;,\;
S_{g_{i_1}^{j_1}g_{i_2}^{j_2}\,...\; g_{i_r}^{j_r}(j)}>$$
is the star $S_{j_2}.$
In particular, the path $\Gamma_{a,b}$ is contained in
$\Omega.$ 
Since $c_1\in V_{M,\infty},$  
the loop $l_1$ does not appear
among the loops $l_{i_1},l_{i_2},\,...\,l_{i_r}.$
Therefore, the common vertex of any two successive stars in the chain
$\Omega$ is not contained in the set $P^{-1}(c_1).$
In particular, among of vertices of $\Gamma_{a,b}$ there are no preimages
of $c_1$ distinct from $a,b$ and hence
$w(1)=1$ on $\Gamma_{a,b}.$

To finish the proof notice that the condition of the
corollary is compositionally stable. Indeed, if 
$L$ is a curve connecting points $a,b$
such that $c_1=P(a)=P(b)$
is a simple point of $P(L)$ and $c_1\in \partial V_{P(L),\infty},$
then $W(L)$ is a curve connecting points $W(a),W(b)$
such that $c_1=\tilde P(W(a))=\tilde P(W(b))$
is a simple point of $\tilde P(W(L))$ and $c_1\in \partial V_{\tilde
P(W(L)),\infty}.$

\vskip 0.2cm

\section{On functional equations
\e{e1} and \e{e2}}

\subsection{Derivation of \e{e1} and \e{e2}}
Let $U$ be a simply connected domain containing no critical values of $P(z)$
such that $S\setminus\{c_1,c_2, ... ,c_{\tilde k}\}\subset U.$
Denote by 
$P_{a_1}^{-1}(z),$
$P_{a_2}^{-1}(z), ... , P_{a_{f_a}}^{-1}(z)$
(resp. $P_{b_1}^{-1}(z),$ $P_{b_2}^{-1}(z),$ ... , 
$P_{b_{f_b}}^{-1}(z)$)
the branches
of $P^{-1}(z)$ in $U$ 
which map points close
to $P(a)$ (resp. $P(b)$)
to points close to $a$ (resp. $b$). 
In particular, $f_a$ (resp. $f_b$)
equals the multiplicity of the point $a$ (resp. $b$) 
with respect to $P(z).$ It was shown in \cite{pa2} for $P(a)=P(b)$
and in \cite{pry} in general case that condition \e{1} implies
equality \e{e1} or system \e{e2}, where
as above $Q(z)=\int q(z) \d z$ is normalized
by the condition ${Q(a)=Q(b)=0}.$

For the sake of self-containedness of this paper below
we provide a short derivation of \e{e1}, \e{e2} from
theorem \ref{t1}.

\bp \l{brc} Suppose that condition \eqref{1} holds. Then, if $P(a)=P(b),$ 
equation \e{e1}
holds in $U$.
Furthermore, if $P(a)\neq P(b),$ then system \e{e2} holds
in $U.$
\ep

\pr Suppose first that $P(a)= P(b)=c_1.$ Examine the relation 
$$\phi_1(z)=\sum_{i=1}^nf_{1,i}Q(P^{-1}_{i}(z))=0.$$ 
Let $i,$ $1\leq i \leq n,$ be an index such that $f_{1,i}\neq 0$
and let $x$ be a vertex of the star $S_i$ such that $P(x)=c_1.$
Observe that if $x\neq a,$ $x\neq b,$ then there exists an index
$\tilde i$ such 
that $x$ is also a vertex of the star $S_{\tilde i}$ and 
$f_{1,\tilde i}=-f_{1,i}.$ Furthermore, we have $\tilde i=g_1^j(i)$
for some natural number $j.$ Therefore, $\phi_1(z)$ has the form
$$\phi_1(z)=-Q(P^{-1}_{i_a}(z))+ $$
$$
Q(P^{-1}_{i_1}(z))-Q(P^{-1}_{g_1^{j_1}(i_1)}(z))+
\, ... \, +
Q(P^{-1}_{i_r}(z))-Q(P^{-1}_{g_1^{j_r}(i_r)}(z))$$
$$+Q(P^{-1}_{i_b}(z))=0,$$
where $i_a$ (resp. $i_b$) is the index such that 
$a\subset S_{i_a}$
(resp. $b\subset S_{i_b}$), $i_1,i_2,\, ... \, i_r$
are some other indices and $j_1,j_2,\, ... \, j_r$
are natural numbers.

Let $n_1$ be the order of the element $g_1$ in the group $G_P.$
For each $s,$ \linebreak ${0\leq s \leq n_1-1,}$ the equality
$$-Q(P^{-1}_{g^s_1(i_a)}(z))+ $$
$$
Q(P^{-1}_{g^s_1(i_1)}(z))-Q(P^{-1}_{g_1^{j_1+s}(i_1)}(z))+
\, ... \, +
Q(P^{-1}_{g^s_1(i_r)}(z))-Q(P^{-1}_{g_1^{j_r+s}(i_r)}(z))$$
$$+Q(P^{-1}_{g^s_1(i_b)}(z))=0$$ holds by the analytic continuation
of the equality 
$\phi_1(z)=0.$ Summing these equalities and taking into account that
for any $i,$ $1\leq i \leq n,$ and any natural number $j$ we have:
$$\sum_{s=0}^{n_1-1} 
Q(P^{-1}_{g^s_1(i)}(z))=\sum_{s=0}^{n_1-1} Q(P^{-1}_{g_1^{j+s}(i)}(z)),$$
we obtain equality \e{e1}.

In the case when $P(a)\neq P(b)$ the proof is similar: if $P(a)=c_1,$
$P(b)=c_2,$ then one must examine relations $\phi_1(z)=0$ and $\phi_2(z)=0.$

\vskip 0.2cm
Note that if points $a,b$ are not critical
points of $P(z),$ then \eqref{e1} reduces to \eqref{eb}
while \eqref{e2} leads to the equality $q(z)\equiv 0.$
In view of lemmas \ref{lcomp} and \ref{rl} this implies
immediately the following result from \cite{c}
(see also \cite{pa2},\cite{pry}).

\bc \l{cris}
Let $P(z),q(z)\in \C[z]$, $q(z)\neq 0,$
$a,b\in \C,$ $a\neq b.$ Suppose that
$a,$ $b$ are not critical points  of $P(z).$
Then conditions \e{1} and \e{2} are 
equivalent. \ec

\subsection{Relations between branches of $Q(P^{-1}(z))$}
In this subsection we examine how linear relations between branches of \linebreak $Q(P^{-1}(z))$ over $\C$ reflect on the structure of coefficients of Puiseux
expansions of $Q(P^{-1}(z))$ near
infinity.

Let $P(z)$ be a non-constant polynomial of degree $n$
and 
let $z_0\in \C$ be a non-critical
value of $P(z).$
If
$\vert z_0\vert $ is sufficiently large then in a neighborhood $U_{z_0}$ 
of $z_0$ each
branch of $P^{-1}(z)$ can be 
represented
by a Puiseux series centered at infinity. More precisely,
if $P_0^{-1}(z)$ is a fixed branch of $P^{-1}(z)$ near $z_0$ then
in $U_{z_0}$ we have: $$P_0^{-1}(z)=
\sum_{k=-1}^{\infty}v_kz^{-\frac{k}{n}}, \ \ \ \ \ \ \ \ v_k\in \C, \ \ \ 
\ \ \ \ \ \varepsilon_n=exp(2\pi i/n),$$
where $z^{\frac{1}{n}}$ is a branch of
the algebraic function which is inverse to $z^n$ in $U_{z_0}.$  
If $l$ is a loop around infinity then 
the result of the analytic continuation of the branch 
$P^{-1}_0(z)$ along $l^j,$ $0\leq j \leq n-1,$ is represented
by the series 
\be \l{cn}
P^{-1}_j(z)=\sum_{k=-1}^{\infty}v_k\varepsilon_n^{jk}z^{-\frac{k}{n}}.
\ee
The numeration of branches of $P^{-1}(z)$ near $z_0$ defined by equation
\eqref{cn} 
is called canonical. Clearly, such a numeration depends on
the choice of $P^{-1}_0(z).$ Nevertheless, any canonical numeration 
induces
the same cyclic ordering of branches of $P^{-1}(z)$ in $U_{z_0}.$
This cyclic ordering also will be called canonical.
For any non-zero polynomial $Q(z),$ $\deg Q(z)=m,$  
the composition $Q(P^{-1}_j(z)),$  $0\leq j \leq n-1,$
is represented near $z_0$ by the series \e{ps2}
obtained by the substitution of series \eqref{cn} in $Q(z).$

Let $U$ be a simply-connected domain containing no critical values of
$P(z)$ such that some linear combination of branches
of $Q(P^{-1}(z))$ over $\C$ identically vanishes in $U$.
Considering in case of necessity a bigger domain we can 
suppose without loss of generality that $\infty \in \partial U.$
Then series \e{cn} converge in a domain
$V\subset U.$
Furthermore, we can assume that
the numeration of branches of $P^{-1}(z)$ in $U$ 
is induced by a canonical numeration
of branches of $P^{-1}(z)$ in $V$.
If equality
\be \l{x}
\sum_{i=0}^{n-1}f_jQ(P^{-1}_j(z))=0, \ \ \ \ \ \ f_j\in \C, 
\ee
holds in $U,$ then substituting in \e{x} expansions \e{ps2}
we see that \e{x} reduces to the system
$$\sum_{j=0}^{n-1}f_ju_k\varepsilon^{kj}_n=0, \ \ \ \ \ \ \ k\geq -m.$$
Introducing the notation $F(z)=\sum_{j=0}^{n-1}f_jz^j$
and summing up we get:

\bl \l{la} The equality \e{x}
holds in $U$ if and only if 
for any $k\geq -m$ 
either
$u_k=0$
or $F(\varepsilon^k_n)=0.$
\el

In particular, since all $u_k$ can not
vanish
and $\deg F(z)<n,$ the following statement is true.

\bc \l{la+}
If equality \e{x} holds in $U$
%for non-constant polynomials $P(z),Q(z),$
then $F(\varepsilon^r_n)=0$ for at least one $r,$ ${0\leq r \leq n-1.}$ On the other hand, for at least one $r,$ $0\leq r \leq n-1,$ the equality $u_k=0$ holds whenever $k\equiv r \ \mod n.$
\ec

\subsection{Lemma about monodromy groups of polynomials}
In order to relate \e{e1}, \e{e2}
with coefficients of Puiseux
expansions of $Q(P^{-1}(z))$ near
infinity
we are going to examine
which roots of unity can be roots of the corresponding polynomial
\be \l{e11}
r(z)=\frac{1}{d_{a}}\sum_{s=1}^{d_{a}}
z^{a_s}-\frac{1}{d_{b}}\sum_{s=1}^{d_{b}}
z^{b_s},
\ee or common roots of the corresponding pair of polynomials
\be \l{e22}
r_1(z)=\frac{1}{d_{a}}\sum_{s=1}^{d_{a}}
z^{a_s}, \ \ \ \ \ \ \ \ \ \ r_2(z)=\frac{1}{d_{b}}\sum_{s=1}^{d_{b}}
z^{b_s}.
\ee
For this propose we establish now a 
geometric property of monodromy groups of polynomials 
which concerns the mutual arrangement
of indices $a_1,$
$a_2, ... , a_{f_a}$
and $b_1,$ $b_2,$ ... , 
$b_{f_b}$ under assumption that the numeration of branches is canonical.

Let $P(z)\in \C[z],$ $\deg P(z)=n,$ $a,b\in \C,$ $a\neq b.$
Let $U$
be a simply-connected domain containing no critical values of $P(z)$ 
such that $P(a),P(b), \infty \in \partial U.$ 
Fix a canonical numeration of branches of $P^{-1}(z)$ in $U$
and let $P_{u_1}^{-1}(z),$
$P_{u_2}^{-1}(z), ... ,$ $P_{u_{f_a}}^{-1}(z)$
(resp. $P_{v_1}^{-1}(z),$ $P_{v_2}^{-1}(z), ... ,
P_{v_{f_b}}^{-1}(z)$) be
the branches of $P^{-1}(z)$ in $U$ which map points close
to $P(a)$ (resp. $P(b)$)
to points close to the point $a$ (resp. $b$) numbered
by means of this numeration.
The lemma below describes the
mutual position on the unit circle of the sets
$V(a)=\{ \v_n^{a_1},\v_n^{a_2}, ... , \v^{a_{f_a}}_n \}$ and
$V(b)=\{\v_n^{b_1},\v_n^{b_2}, ... ,\v_n^{b_{f_b}}\},$
where $\varepsilon_n=exp(2\pi i/n).$
Let us introduce the following
definitions.
Say that two sets of points $X,Y$ on the unit circle $S_1$ are
{\it disjointed} if there exist $s_1, s_2 \in S_1$
such that all points from $X$ are on
the one of two connected components of $S_1\setminus \{s_1,
s_2\}$
while all points from $Y$ are on
the other one. Say that $X,Y$ are
{\it almost disjointed} if $X\cap Y$ consists of a single point $s_1$
and there exists a point $s_2\in S_1$ such that
all points from $X\setminus s_1$ are on
the one of two connected components of $S_1\setminus \{s_1,
s_2\}$ while all points from $Y\setminus s_1$ are on
the other one.

\vskip 0.2cm

\noindent{\bf Monodromy Lemma.} {\it 
The sets $V(a)$ and $V(b)$ are almost disjointed. Furthermore,
if $P(a)= P(b)$ then $V(a)$ and $V(b)$ are disjointed.
}

\vskip 0.2cm
 
\noindent{\it Proof.}
Consider first the case when $P(a)=P(b).$
\begin{figure}
\medskip
\epsfxsize=9truecm
\centerline{\epsffile{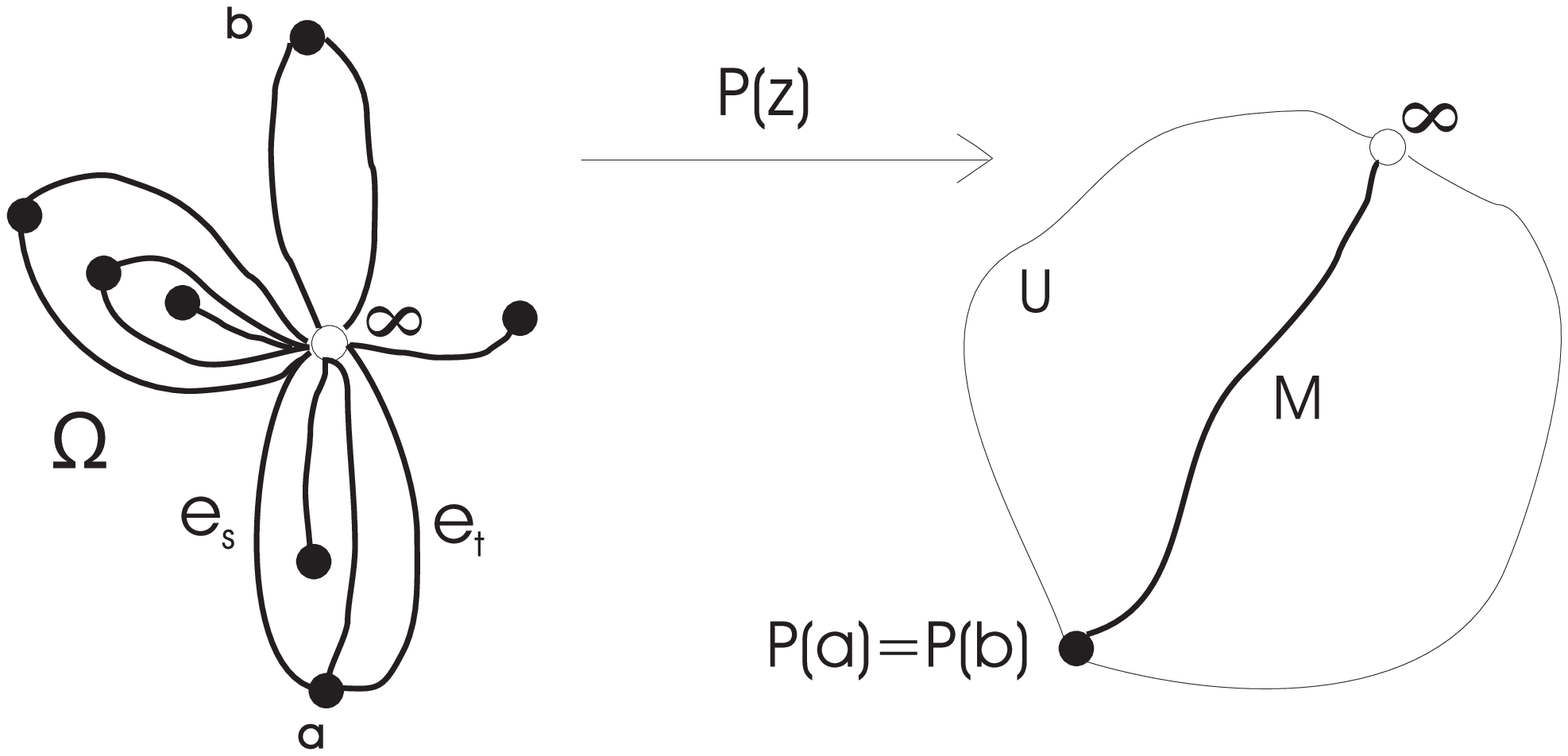}}
\smallskip
\centerline{Figure 3}
\medskip
\end{figure} 
Let $M\subset U$ be a simple curve connecting points $P(a)=P(b)$
and $\infty.$
Consider the preimage
$P^{-1}\{M\}$ of $M$ under the map $P(z)\,:\,\C\P^1\rightarrow
\C\P^1.$
It is convenient to consider $P^{-1}\{M\}$ as a bicolored graph
$\Omega$ embedded into the Riemann sphere:
the black vertices of $\Omega$ are preimages of
$P(a)=P(b),$ the
unique white vertex is the preimage of $\infty,$ and
the edges of $\Omega$ are preimages of $M$ (see Fig. 3). 
Since the multiplicity of the vertex $\infty$ equals $n$ and $\Omega$
has $n$ edges, $\Omega$ is connected.
The edges of $\Omega$ are identified with branches of $P^{-1}(z)$ in $U$
as follows:
to a branch $P^{-1}_k(z),$ $1\leq k \leq n,$  
corresponds the edge $e_k$
such that $P^{-1}_k(z)$ maps $M\setminus \{P(a),\infty\}$ into $e_k.$
For any vertex $v$ of $\Omega$ the orientation of $\C\P^1$
induces a natural cyclic ordering on edges of $\Omega$ adjacent to
$v.$ In particular, taking $v=\infty,$ we obtain a 
cyclic ordering on edges of $\Omega.$
Clearly, this cyclic ordering coincides with that induced by
the canonical cyclic ordering of branches of $P^{-1}(z)$ in $U.$
Let $E_a=\{e_{a_1},e_{a_2},\, ... \,, e_{a_{f_a}}\}$
(resp. $E_b=\{e_{b_1},e_{b_2},\, ... \,, e_{b_{f_b}}\}$) be the union
of edges of
$\Omega$ which are adjacent to the vertex $a$ (resp. $b$).
Let $D$
be the domain from the collection of domains
$\C\P^{1}\setminus E_a$ which contains point $b$
and let $e_s,e_t\in E_a$ be the edges which bound $D.$
Clearly, all the edges from $E_a$ are contained in
$\C\P^1\setminus D.$
Therefore, the lemma is equivalent to the following statement: 
the domain $D$ contains $e_h\setminus {\infty}$
for all $e_h\in E_b.$ But the last statement is a corollary
of the Jordan theorem since an edge $e_h\in E_b$ can intersect $e_s$ or   
$e_t$ only at infinity.

In the case when $P(a)\neq P(b)$ the proof is modified as follows.
Divide the boundary of $U$ into three
parts $M_1,M_2,M_3,$ where $M_1$ connects the point $\infty$ with the 
point 
$P(a),$ $M_2$ connects the point $\infty$ with the point 
$P(b),$ and $M_3$ connects the point $P(a)$ with the point 
$P(b).$ 
Consider now $P^{-1}\{\partial U\}$ as a graph
$\Omega$ embedded into the Riemann sphere. The vertices of $\Omega$ 
are divided into three groups: the first one consists of vertices  
that are preimages of $\infty,$ the second one consists of vertices  
that are preimages of $P(b),$ and the third one consists of vertices  
that are preimages of $P(a).$ Similarly, the edges of $\Omega$ 
also are divided into three groups: the first one consists of edges
that are preimages of $M_1,$ the second one consists of edges 
that are preimages of $M_2,$ and the third one consists of edges 
that are preimages of $M_3.$ Finally, the faces of $\Omega$ are divided 
into two groups: the first one consists of faces
that are preimages of $U$ and the second one consists of faces
that are preimages of $\C\P^1\setminus \overline{U}$
(see Fig. 4).

\begin{figure}
\medskip
\epsfxsize=11truecm
\centerline{\epsffile{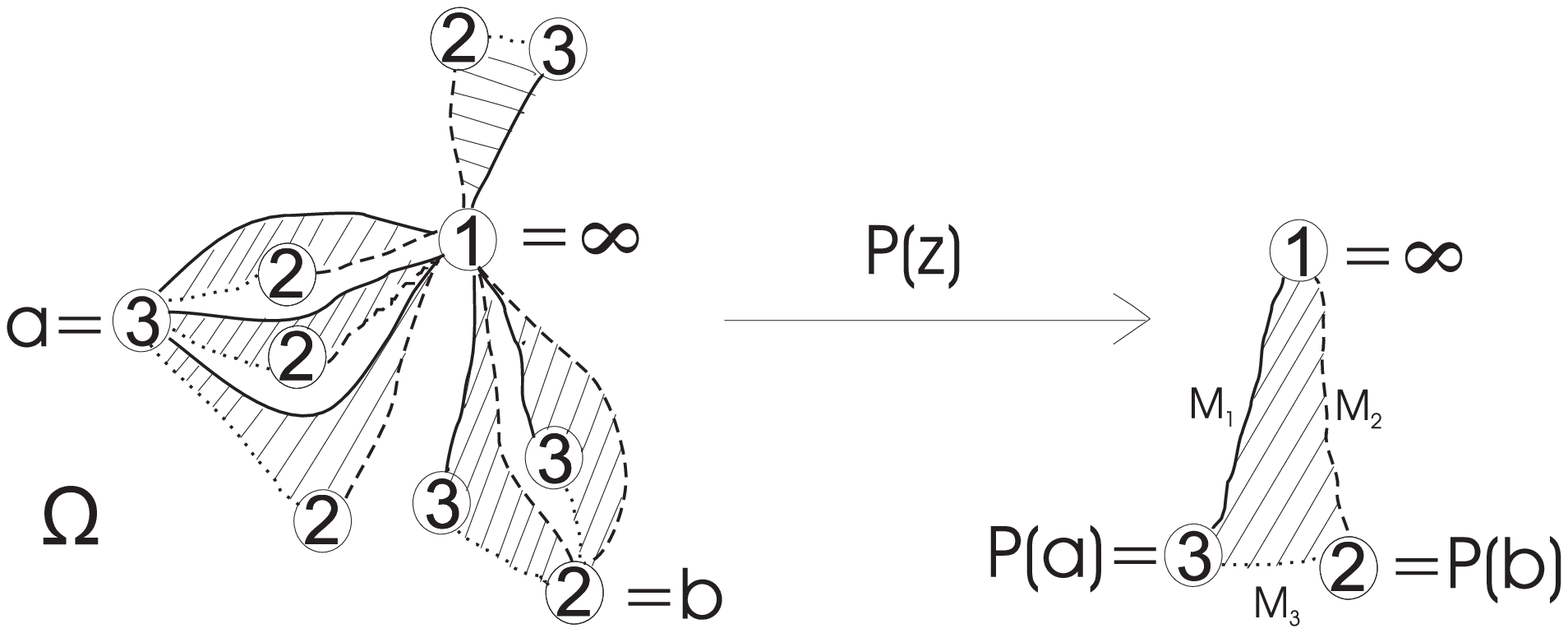}}
\smallskip
\centerline{Figure 4}
\medskip
\end{figure}

The faces from the first group are identified with branches of $P^{-1}(z)$ 
in $U$
as follows:
%if $N_3$ is a numeration of branches of $P^{-1}(z)$ in $V$
to a branch $P^{-1}_k(z),$ $1\leq k \leq n,$  
corresponds the face $f_k$
such that $P^{-1}_k(z)$ maps bijectively $U$ on $f_k\setminus \partial f_k.$
The edges from the corresponding groups which bound $f_k$ will be denoted 
by
$e^1_k, e^2_k, e^3_k$ correspondingly.  
Note that in a counterclockwise direction around infinity the edge 
$e^1_k,$ $1\leq k \leq n,$ is
followed by the edge $e^2_k.$
The canonical cyclic ordering of
branches of $P^{-1}(z)$ in $U$ induces a cyclic ordering of faces of the
first group of $\Omega$ and this  
ordering coincides with that induced by the
orientation of $\C\P^1.$
Let $E_a^1=\{e^1_{a_1},e^1_{a_2},\, ... \,, e^1_{a_{f_a}}\}$
(resp. $E_b^2=\{e^2_{b_1},e^2_{b_2},\, ... \,, e^2_{b_{f_b}}\}$)
be the union of edges from the first (resp. 
the second) group
$\Omega$ which are adjacent to the vertex $a$ (resp. $b$).
Let $D$
be the domain from the collection of domains
$\C\P^{1}\setminus E_a^1$ which contains point $b.$
%and let $e_s^1,e_t^1\in E_a$ be the edges which bound $D.$
Once again the Jordan theorem implies that
all the edges from $E_a^1$ are contained in
$\C\P^1\setminus D$ while
$D$ contains $e_h^2\setminus {\infty}$
for all $e_h^2\in E_b^2.$ Taking into account that
for any $k,$ $1\leq k \leq n,$
the edge $e^1_k$ is 
followed by $e^2_k$ this fact implies
that $V(a)$ and $V(b)$ are almost disjointed.
Note that, in contrast to the case when 
$P(a)=P(b),$
now the sets $V(a)$ and $V(b)$ can have a non-empty intersection
consisting of a single element.

\subsection{On coefficients of Puiseux expansion of $Q(P^{-1}(z))$}

In this subsection we deduce from the monodromy lemma the following
important
property of
the Puiseux expansions \e{ps2}
for pairs $P(z),Q(z)$
satisfying \e{e1}, \e{e2}.

\bt \l{pu}
Let $P(z),Q(z)\in \C[z]$, $\deg P(z)=n,$ 
$a,b\in \C,$ $a\neq b.$
Suppose that \e{e1} or \e{e2} holds.
Then $u_k=0$ for any $k$ such that $\GCD(k,n)=1.$ 

\et

\pr 
Suppose first that $P(a)=P(b).$ 
Then lemma \ref{la} implies that 
$u_k = 0$ whenever the number $\varepsilon_n^{k}$ is not
a root of the polynomial \e{e11}.
Let us show that if $\GCD(k,n)=1$ 
then
the equality $r(\varepsilon_n^{k})=0$ is impossible.
Indeed, if $(k,n)=1,$
then $\varepsilon_n^{k}$
is a primitive $n$-th root of unity.
Since
the $n$-th cyclotomic polynomial $\Phi_n(z)$ is irreducible over $\Z$,
the equality $r(\varepsilon_n^{k})=0$ 
implies that $\Phi_n(z)$ divides $r(z)$ in the ring $\Z[z].$
Therefore, the primitive $n$-th root of unity
$\varepsilon_n=exp(2\pi i/n)$ also is
a root of $r(z)$ and hence the equality   
$$\sum_{s=1}^{d_a}\varepsilon_n^{a_s}/d_a =
\sum_{s=1}^{d_b}\varepsilon_n^{b_s}/d_b$$ holds.
The last equality is equivalent to the statement
that the mass centers of the sets $V(a)$ and $V(b)$ 
coincide. But this contradicts to the monodromy lemma. Indeed, 
the mass center of a system of points in $\C$
is inside of the convex envelope of this system and therefore
the mass centers of disjointed sets must be distinct.

If $P(a)\neq P(b)$ then, similarly, the inequality $u_k\neq 0$ for
$\GCD(k,n)=1$ implies that
$$\sum_{s=1}^{d_a}\varepsilon_n^{a_s}/d_a =0, \ \ \ \ \ 
\sum_{s=1}^{d_b}\varepsilon_n^{b_s}/d_b=0.$$
But this again contradicts the monodromy lemma. Indeed, the fact that  
the sets $V(a)$ and $V(b)$ are almost disjoint implies that at least 
one
from these sets is contained in an open half plane bounded by a line
passing through the origin and therefore the mass center of this set
is distinct from zero.

\bc
Let $P(z),q(z)\in \C[z]$, $q(z)\neq 0,$ $\deg P(z)=n,$ 
$\deg Q(z)=m,$ 
$a,b\in \C,$ $a\neq b.$ Suppose that \e{1} holds.
Then $\GCD(m,n)>1.$
\ec

\pr Since in expansions \e{cn}
the coefficient $v_{-1}$ is distinct from zero,
the coefficient $u_{-m}=v_{-1}^m$ in expansions \e{ps2}
is distinct from zero.
Since \e{1} implies \e{e1} or \e{e2}
by proposition \ref{brc} it follows from
theorem \ref{pu} that $\GCD(m,n)>1.$

\vskip 0.2cm
Notice that theorem \ref{pu} agrees with conjecture
\e{cc}. Indeed, if
\be \l{a}
Q(z)=\tilde Q_1(W_1(z))+\tilde Q_1(W_1(z))+ \ ... \ + \tilde Q_r(W_r(z)),
\ee
where $W_1(z),W_2(z),...,W_r(z)$ are (non-trivial) right divisors of
$P(z)$ in the composition algebra, 
$$P(z)=\tilde P_1(W_1(z))=\tilde P_2(W_2(z))=\ ... \ =
\tilde P_r(W_r(z)),$$
%for non-linear
%$W_1(z),W_2(z),...,W_r(z)$$W_1(z),W_2(z),...,W_r(z)$
then
the expansion \e{ps2} has the form
$$Q(P^{-1}(z))=\tilde Q_1(\tilde P_1^{-1}(z))+\tilde Q_2(\tilde
P_2^{-1}(z))+\
... \ +\tilde Q_r(\tilde P_r^{-1}(z)).$$
Since $\deg \tilde P_j(z) < n,$ $1\leq j \leq r,$ it follows
easily that
$u_k=0$ for any $k$ such that $\GCD(k,n)=1.$
Conjecturally, vice versa, equalities
$u_k=0$ for all $k$ with $\GCD(k,n)=1$ imply
that $Q(z)$ has form \e{a} at least under some additional
assumptions. We plan to discuss this topic in another paper.

\section{Further description of definite polynomials}

\subsection{Case when $a$ or $b$ is not a critical
point of $P(z)$}
As a first application of the Puiseux expansions technique 
we provide in this section the following generalization of
corollary \ref{cris}.

\bt \l{alo} Let $P(z),q(z)\in \C[z]$, $q(z)\neq 0,$ 
$a,b\in \C,$ $a\neq b.$ Suppose that at least one from 
points $a$ and $b$ is not a critical point of the polynomial $P(z).$
Then
conditions \e{1} and \e{2} are equivalent.
\et 
\pr Since the condition of the theorem is compositionally
stable it follows from lemmas \ref{rl}, \ref{lcomp}
that we only must show that equality \e{eb} holds.
To be definite suppose that the point $a$ is not a critical point of 
$P(z).$
By proposition \ref{brc} either the system
\be \l{pk}
Q(p_{a_1}^{-1}(z))=0, \ \ \ \ \ \ \ \  \sum_{s=1}^{d_{b}}
Q(p_{b_s}^{-1}(z))=0
\ee
or the equality
\be \l{pk2}
Q(p_{a_1}^{-1}(z))=\frac{1}{d_{b}}\sum_{s=1}^{d_{b}}
Q(p_{b_s}^{-1}(z))
\ee
holds. Nevertheless, since the first equation of system \e{pk}
leads to the equality ${q(z)\equiv 0,}$ we must only consider
equation \e{pk2}.

Applying lemma \ref{la} we see 
that for any $k$ such that $u_k\neq 0$
the equality  
$$d_b(\varepsilon_n^k)^{a_1}=\sum_{s=1}^{d_{b}} (\varepsilon_n^k)^{b_s} $$ 
holds. 
The triangle inequality implies that
this is possible only if 
$$(\varepsilon_n^k)^{a_1} =(\varepsilon_n^k)^{b_1}=(\varepsilon_n^k)^{b_2}= 
... =(\varepsilon_n^k)^{b_{d_b}}.$$ 
Therefore,
$$Q(p_{a_1}^{-1}(z))=
Q(p_{b_1}^{-1}(z))=Q(p_{b_2}^{-1}(z))= ... =
Q(p_{b_{d_s}}^{-1}(z))
.$$

\subsection{Case when
$\deg P(z)=p^n$}
In this subsection we deduce from theorem \ref{pu} the solution
of the polynomial moment problem in the case
when $\deg P(z)=p^r$ for $p$ prime.

\bt \l{prime} Let $P(z),q(z)\in \C[z],$ $q(z)\neq 0,$ 
$a,b\in \C,$ $a\neq b.$
Suppose that $\deg P(z)=p^r,$
where $p$ is a prime number. Then conditions \e{1} and \e{2} are equivalent.
\et

\pr Again, since the condition of the theorem is compositionally
stable, it is enough to show that \e{eb} holds. Consider
expansion \e{ps2}. By theorem \ref{pu} the equality
$u_k=0$ holds for any $k$ with $\GCD(k,p^r)=1.$
Show that this fact implies the equality
$$
Q(P^{-1}_j(z))=Q(P^{-1}_{j+p^{r-1}}(z))
$$
for any $j,$ $0\leq j\leq n-1.$
Indeed, we have:
$$ 
Q(P^{-1}_j(z))-Q(P^{-1}_{j+p^{r-1}}(z))
=\sum_{k=-m}^{\infty}
w_kz^{-\frac{k}{n}},
$$
where
$$
w_k=u_k(\varepsilon_{p^r}^{jk}-\varepsilon_{p^r}^{(j+p^{r-1})k}).
$$
If $\GCD(k,p^r)=1$ then $u_k=0$ and hence $w_k=0.$
Otherwise, $k=p\tilde k$ for some $\tilde k\in \Z.$ 
Therefore,
$$\varepsilon_{p^r}^{(j+p^{r-1})k}=
\varepsilon_{p^r}^{jk}\varepsilon_{p^r}^{p^{r}\tilde k}=
\varepsilon_{p^r}^{jk}$$ and hence again
$w_k=0.$

\subsection{Case when $P(z)$ is indecomposable}
%In this subsection we give, using the obtained results,
Theorems \ref{t1}, \ref{prime} allow us to give 
a short proof of the theorem proved in \cite{pa2}, \cite{pa4}
which describes solutions to \e{1} in case when $P(z)$
is indecomposable that is
can not be represented as a composition
$P(z)=P_1(P_2(z))$ with non-linear polynomials $P_1(z),$
$P_2(z).$
\bt \l{nepr}
Let $P(z),q(z)\in \C[z],$ $q(z)\neq 0,$ 
$a,b\in \C,$ $a\neq b.$
Suppose that $P(z)$ is indecomposable.
Then conditions \e{1} and \e{2} are equivalent. In more details,
$Q(z)$ is a polynomial in $P(z)$ 
and $P(a)=P(b).$
\et

\pr 
Once again we only must prove that \e{eb} holds.
Suppose the contrary that is that
all $Q(P^{-1}_i(z)),$ $1\leq i \leq n,$ where $n=\deg P(z)$ are
different; then the monodromy group $G$ of the algebraic
function $Q(P^{-1}(z))$ obtained by the complete analytic continuation
of $Q(P^{-1}_i(z)),$ $1\leq i \leq n,$ 
coincides with that of $P^{-1}(z)$.
Since $P(z)$ is indecomposable, $G$ is primitive
by the Ritt theorem \cite{ri}.
Since for the case when $n=\deg P(z)$ is a prime number the statement follows
from theorem \ref{prime} we can suppose that $n$ is a composite number.
By the Schur theorem (see e.g. \cite{wi}, Th. 25.3)
a primitive permutation group of composite degree $n$
which contains an $n$-cycle is doubly transitive.

Recall now the following fact:
roots
$\alpha_i,$ $1\leq i \leq n,$ of an
irreducible algebraic equation over a field $k$ of characteristic zero
with doubly transitive Galois group can not satisfy any
relation $$\sum_{i=1}^{n} c_i\alpha_i=0,\ \ \ \ \ \ c_i\in k,$$ except
the case when $c_1=c_2= ... =c_n$ 
(see \cite{gi}, Proposition 4, or, in the context of algebraic functions, \cite{pa2}, lemma 2).
Since the monodromy group of an algebraic function coincides with
the Galois group of the equation over $\C(z)$
which defines this function, 
it follows that if all $Q(P^{-1}_i(z)),$ $1\leq i \leq n,$ are
different, then equality \e{x} is possible only when
\be \l{ko}
f_1=f_2= ... =f_n.
\ee 
On the other hand, for any non-trivial equation $\phi_s(z)=0$
appeared in theorem \ref{t1}
the equality \e{ko} is impossible by construction.
This contradiction completes the proof.

\section{Solution of the polynomial moment problem
for polynomials of degree less than 
\protect
$10$}

In this section
%using obtained results combined with some combinatorial
%reasoning involving cacti
we provide a complete solution of 
the polynomial moment problem
for polynomials of degree less than 
$10.$

For an extended cactus $\tilde \lambda_P$ and a path $\Gamma_{a,b}$ define
{\it the 
skeleton}
$\hat\Gamma_{a,b}$ of $\Gamma_{a,b}$
as follows.
Draw the path
$\Gamma_{a,b}$ separately from the graph $\tilde\lambda_P$ 
and erase all its white vertices.
Number the edges of the obtained graph $\hat\Gamma_{a,b}$ so that
the number of an edge $e_k$ coincides with
the number of the star $S_k$ of $\tilde\lambda_P$
for which $e_k\subset S_k.$   
The number of edges of
$\hat\Gamma_{a,b}$ is called the length
$l(\hat\Gamma_{a,b})$ of $\hat\Gamma_{a,b}.$
For example, the skeleton 
$\hat\Gamma_{a,b}$ of the path $\Gamma_{a,b}$ from Fig. 1
is shown on Fig. 5; here $l(\hat\Gamma_{a,b})=4.$

\medskip
\epsfxsize=8truecm
\centerline{\epsffile{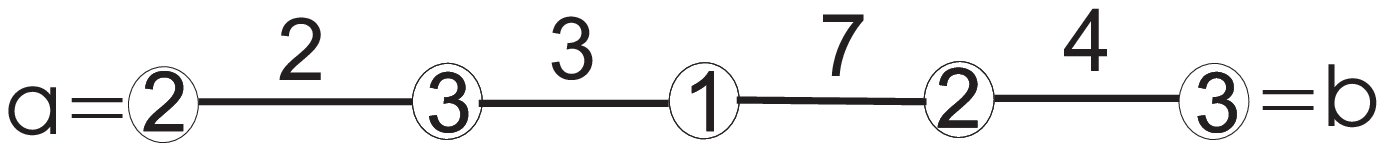}}
\smallskip
\centerline{Figure 5}
\medskip

\bt 
Let $P(z),q(z)\in \C[z]$, $q(z)\neq 0,$ 
$a,b\in \C,$ $a\neq b,$ satisfy \e{1}.
Suppose that $\deg P(z)<10.$ Then
either condition \e{2} holds  or
%$\deg P(z)=6$ and
there exist linear
functions $L_1(z),L_2(z)$ such that
$$L_2(P(L_1(z)))=T_6(z), \ \ \ \ L_1^{-1}(a)=-\sqrt{3}/2,
\ \ \ \ L_1^{-1}(b)=\sqrt{3}/2,$$
and $$Q(L_1(z))=A(T_3(z))+B(T_2(z))$$
for some $A(z),B(z)\in \C[z].$
\et

\pr First of all observe that any natural number $n<10$ distinct
from 6
is either a prime number or a degree of a prime number.
Therefore, it follows from theorem \ref{prime}
that
it suffices to consider the case
when $\deg P(z) = 6.$ Furthermore, in view of theorem \ref{alo}
we can suppose that the points $a,b$ are critical points of $P(z).$
Finally notice that in order to prove that condition \e{2}
holds for $P(z),$ $q(z)$ satisfying \e{1} with $\deg P(z)=6$
it is enough to establish equality \e{comp}.
Indeed, if $W(a)\neq W(b)$ in \e{comp} then performing the change of
variable $z\rightarrow W(z)$ we see
that \e{1} holds
for $\tilde P(z), \tilde Q(z), W(a), W(b).$ 
If $\deg W(z)$ equals 3 or 2,
then it follows from theorem 
\ref{nepr} that ${\tilde Q(z)=R(\tilde P(z))}$ for some $R(z)\in \C[z]$ and 
$\tilde P(W(a))=\tilde P(W(b)).$ Therefore, \e{2} holds
with $W(z)=P(z), \tilde Q(z)=R(z).$ On the other hand, if 
$\deg W(z)=6$ in \e{comp} then necessary   
$W(a)= W(b)$ since otherwise $\tilde Q(z)$ would be orthogonal to
all powers of $z$ on the segment $W(a),W(b).$ 
In particular, in view of lemma \ref{lcomp}, we see that
in order to prove that conditions \e{1} and \e{2} are equivalent
it is enough to establish \e{eb}.

Since $\deg P(z)= 6,$ clearly $l(\hat\Gamma_{a,b})\leq 6.$ Moreover,
since the points $a,b$ are critical points of $P(z),$
the valency of the corresponding vertices of $\tilde \lambda_P$ is
at least 2, and, therefore, actually $l(\hat\Gamma_{a,b})\leq 4.$
Consider all possible cases. First of all observe that
the equality $l(\hat\Gamma_{a,b})=1$ is
impossible. Indeed,
in this case theorem \ref{t1} implies that
$Q(P^{-1}_i(z))=0,$ where $i$ is 
the number of the unique edge of $\hat\Gamma_{a,b},$ 
and therefore $q(z)\equiv 0.$
Furthermore, if $l(\hat\Gamma_{a,b})=2$ 
then,
since adjacent vertices of $\hat\Gamma_{a,b}$ have different 
colors,
$\hat\Gamma_{a,b}$ can be of
one from the following two forms shown on Fig. 6.
\vskip 0.2cm
\medskip \epsfxsize=12truecm \centerline{\epsffile{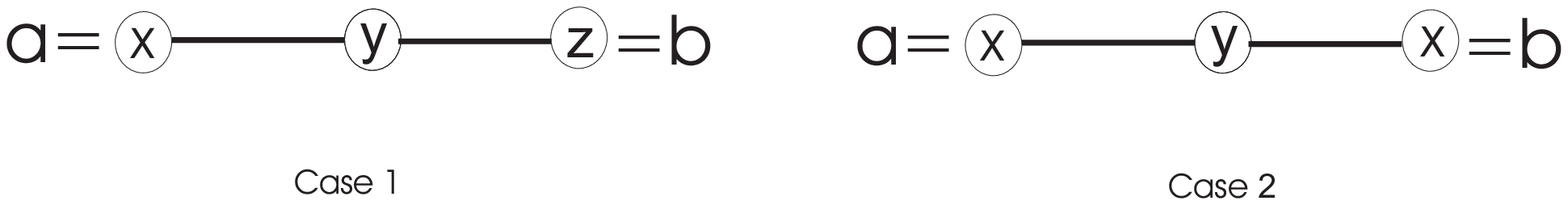}} \smallskip
\centerline{Figure 6} \medskip \noindent
In both cases for the middle
vertex $y$ we have $w(y)=1.$
Therefore by theorem \ref{rt} equality \e{comp} holds and
hence conditions \e{1} and \e{2} are equivalent.
Observe, however, that the first configuration shown on Fig. 6 is actually
not realizable since \e{2} implies that $P(a)=P(b).$

Consider now the case when $l(\hat\Gamma_{a,b})=3.$ It is not
difficult to see that in this case either again $w(y)=1$ for some color
$y$ or $\hat\Gamma_{a,b}$ has the form shown on Fig. 7.

\medskip
\epsfxsize=8truecm
\centerline{\epsffile{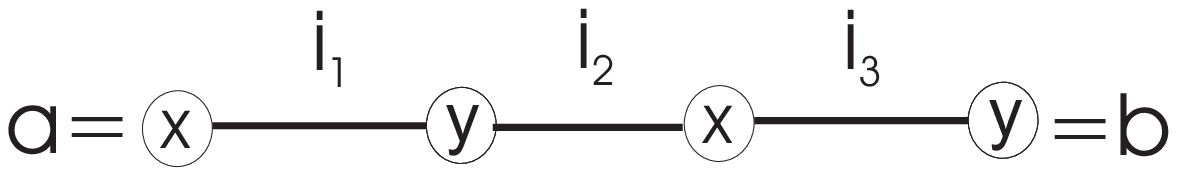}}
\smallskip
\centerline{Figure 7}
\medskip
\noindent
Let us examine the last case.
Since for the skeleton shown on Fig. 7 
we have $P(a)\neq P(b),$ it follows
from proposition \ref{brc} that
system \ref{e2} holds.
Furthermore, the equality $\deg P(z)= 6$ implies that for at least one
point $s,$ $s\in \{a,b\},$
the following two conditions are satisfied: the multiplicity of the 
vertex $s$ of the graph $\tilde \lambda_P$
equals 2 and the connectivity component of $\tilde \lambda_P\setminus s$
which does not contain $\Gamma_{a,b}$ consists of a unique star.
To be definite suppose that $s=a.$
Then, in notation of \ref{cacti},
the first condition implies that \be \l{gra}
\sum_{s=1}^{d_{a}}
Q(P_{a_s}^{-1}(z))=Q(P^{-1}_{i_1}(z))+Q(P^{-1}_{g_x(i_1)}(z))=0 \ee and 
the second one that
$g_y(g_x(i_1))=g_x(i_1).$ Therefore, the analytic continuation of
\eqref{gra} along the loop $l_y$ leads to the equality \be \l{gra2}
Q(P^{-1}_{i_2}(z))+Q(P^{-1}_{g_x(i_1)}(z))=0. \ee Now equalities
\eqref{gra},\eqref{gra2} imply that 
$Q(P^{-1}_{i_1}(z))=Q(P^{-1}_{i_2}(z))$
%Since this equality implies \e{zae}
and we conclude as above that the configuration
shown on Fig. 7 is not realizable.

Consider finally the case when $l(\hat\Gamma_{a,b})=4.$
Since $\hat\Gamma_{a,b}$ has 5 vertices, 
either $w(y)=1$ for some color
$y$ or $\hat\Gamma_{a,b}$ is two-colored. In the last case
$\hat\Gamma_{a,b}$ has the form shown on Fig. 8
\vskip 0.2cm
\medskip
\epsfxsize=8truecm
\centerline{\epsffile{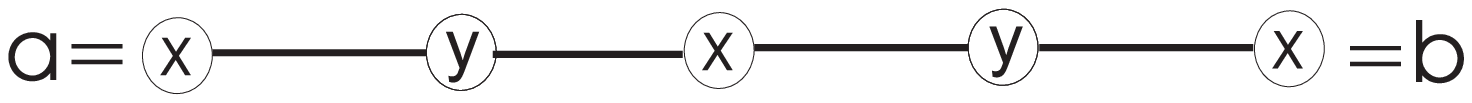}}
\smallskip
\centerline{Figure 8}
\medskip
\noindent
and the corresponding cactus $\tilde \lambda_P$ is a 6-chain
(the cactus with 6 stars of the
maximal diameter).
Furthermore, since $\deg P(z)=6,$ it follows from
the Riemann-Hurwitz formula that
$$\sum_{z\in \C\P^1}({\rm mult}_zP-1)=10.$$ Since
${\rm mult}_{\infty}P-1=5$ and
the combinatorics of $\tilde \lambda_P$ imply that 
$$
\sum_{P(z)=c_x}({\rm mult}_zP-1)=3, \ \ \ \ \ \ \ \
\sum_{P(z)=c_y}({\rm mult}_zP-1)=2,$$
we conclude that $P(z)$ has only two finite critical values $c_x,c_y.$

It follows from the Riemann existence theorem (see e.g. \cite{cac2})
that a complex polynomial with given
critical values 
%$c_x,c_y$ 
is defined by its cactus up to a 
linear change of variable.
On the other hand, it is easy to see using the
formula $T_n(\cos \phi)=\cos n\phi$ that 
$T_n(z)$ has only two critical values $-1,1$ and that
all critical points of $T_n(z)$ are simple,
Therefore, the corresponding cactus is a chain.
In particular, for $P(z)=T_6(z)$ the corresponding cactus realized
as the preimage
of the segment $[-1,1]$ (considered as a star connecting
$0$ with points $1$ and $-1$) has the form shown on Fig. 9
(white vertices are omitted).

%\begin{figure}
\vskip 0.5cm
\medskip
\epsfxsize=9.5truecm
\centerline{\epsffile{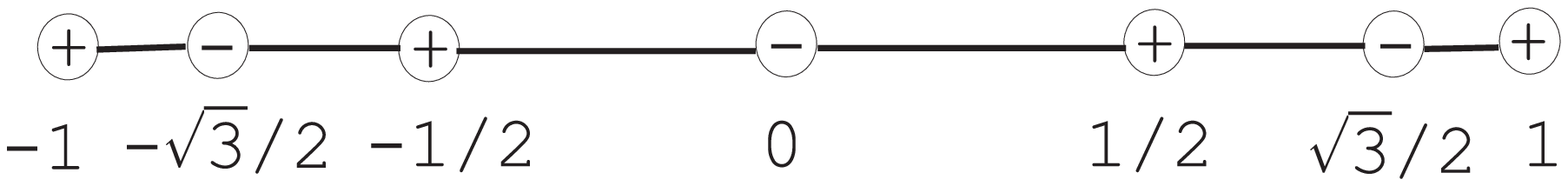}}
%\smallskip
\centerline{Figure 9}
\medskip
\noindent
%\end{figure}
Therefore, if we choose linear functions
$L_1(z), L_2(z)$ such that:
$$L_1^{-1}(a)=-\sqrt{3}/2, \ \ \ L_1^{-1}(b)=\sqrt{3}/2,
\ \ \ 
%critical values of the polynomial
%$L_2(P(L_1(z)))$ are $-1,1,$ and
%$L_2(P(-\sqrt{3}/2))=L_2(P(\sqrt{3}/2))=-1,$
L_2(c_x)=-1,\ \ \ L_2(c_y)=1,$$
the polynomial $L_2(P(L_1(z)))$ will be equal $T_6(z).$

Finally, the last assertion of the theorem follows from
the main result of the paper \cite{pa3} where 
all solutions to \eqref{1} for $P(z)=T_n(z)$
were described.

\bibliographystyle{amsplain}

\end{document}